\documentclass[11pt]{article}

\usepackage{amsfonts}
\usepackage{graphicx}
\usepackage{epstopdf}
\usepackage{amsopn}
\usepackage{amstext,amssymb,amsmath,amsbsy}
\usepackage{amscd}
\usepackage{amsfonts}
\usepackage{indentfirst}
\usepackage{verbatim}
\usepackage{amsmath}
\usepackage{enumerate}
\usepackage{graphicx}
\usepackage{color}
\usepackage[OT1]{fontenc}
\usepackage[latin1]{inputenc}
\usepackage[english]{babel}
\usepackage{amssymb}
\usepackage{subfig}
\usepackage{algorithm}
\usepackage{algpseudocode}
\usepackage{mathtools}
\DeclarePairedDelimiter\ceil{\lceil}{\rceil}

\newcommand{\R}{\mathbb{R}}

\newcommand{\x}{{\bf x}}

\newcommand{\ii}{{\rm i}}
\newcommand{\ik}{\ii k}

\newcommand*{\QEDB}{\null\nobreak\hfill\ensuremath{\square}}
\newtheorem{theorem}{Theorem}[section]
\newtheorem{Remark}{Remark}
\newtheorem{Assumption}{Assumption}

\newtheorem{problem}{Problem}[section]

\numberwithin{equation}{section}

\begin{document}
%%%%% title : short title may not be used but TITLE is required.
% \title{TITLE}
%\title[short title]{TITLE}
\title{The gradient descent method for the convexification to solve boundary value problems of  quasi-linear PDEs  and  a coefficient inverse problem}

 \author{Thuy T. Le\thanks{
Department of Mathematics and Statistics, University of North Carolina at
Charlotte, Charlotte, NC 28223, USA, tle55@uncc.edu, loc.nguyen@uncc.edu
(corresponding author)} \and Loc H. Nguyen\footnotemark[1]}

\date{}
\maketitle

\begin{abstract}
	We study the global convergence of the gradient descent method of the
minimization of strictly convex functionals on an open and bounded set of a
Hilbert space. Such results are unknown for this type of sets, unlike the
case of the entire Hilbert space. 
Then, we use our result to establish a
general framework to numerically solve boundary value problems for
quasi-linear partial differential equations (PDEs) with noisy Cauchy data.
The procedure involves the use of Carleman weight functions to convexify a
cost functional arising from the given boundary value problem and thus to
ensure the convergence of the gradient descent method above. We prove the
global convergence of the method as the noise tends to 0. The convergence
rate is Lipschitz. Next, we apply this method to solve a highly nonlinear
and severely ill-posed coefficient inverse problem, which is the so-called
back scattering inverse problem. This problem has many real-world
applications. Numerical examples are presented.
\end{abstract}



\section{Introduction}

\label{sec Intro}

Numerical solutions of ill-posed Cauchy problems for quasi-linear partial
differential equations (PDEs) is an important topic that arises in many
real-world applications. For example, in the case of parabolic PDEs such
problems are common in heat conduction \cite{Alifanov:sp1994,Alifanov:bh1995}. A natural approach to solve such a problem is to minimize the functional
defined by the least-squares method. However, due to the presence of the
nonlinearity, this functional is non convex. It might have multiple local
minima and ravines. Therefore, a good initial guess, which is located
sufficiently close to the true solution, plays an important role in the
minimization process. 
Since such a good initial guess is not always available,
we, in this paper, use the convexification method, in which it is not
necessary to have a small distance between the starting point of iterations
and the true solution. 
The main content of the convexification method is to
construct a weighted cost functional, which is strictly convex on an a
priori chosen bounded set. It is important that smallness condition is not
imposed on the diameter of this set. 
The unique minimizer  of
that functional on that set is close to the true solution of the given ill-posed Cauchy
problem. The key element of that functional is the Carleman Weight Function
(CWF), which is involved as weight in the Carleman
estimate for the corresponding PDE operator.

An important question arises immediately on how to efficiently find the
global minimizer of such a convex functional on that bounded set. It is well
known that if a functional is strictly convex on the whole Hilbert space,
then the gradient descent method converges to its unique minimizer if
starting from an arbitrary point of that space. 
However, it is not clear
what to do in our case when the strict convexity takes place only on a
bounded set. To address this question, it was proposed in \cite%
{KlibanovNik:ra2017} to use the gradient projection method. However, this
method is a complicated one and is hard to implement numerically. On the
other hand, it was heuristically observed in all numerical studies of the
convexification conducted so far that a simpler gradient descent method
 works well; see e.g., \cite{KlibanovNik:ra2017,VoKlibanovNguyen:IP2020,Khoaelal:IPSE2021,KhoaKlibanovLoc:SIAMImaging2020,KlibanovLiBook}.
This motivates us to analytically study  the question of the global
convergence of the gradient descent method.

More precisely, we prove that the gradient descent method delivers a
sequence converging to the minimizer of that functional on that bounded set
if starting from an arbitrary point of that set. Since smallness conditions
are not imposed on the diameter of this set,  this is global
convergence, see, e.g. \cite{KlibanovNik:ra2017} where the notion of global
convergence is defined. Some numerical results by gradient descent method
will be presented.

Another important part of this paper is to
apply this result to solve a highly nonlinear and severely ill-posed
coefficient inverse problem with a single measurement of back scattering
data in the frequency domain.

As mentioned above, the main idea of the convexification method is to
construct a strictly convex functional. To do this, one uses the Carleman
weight function to convexify the mismatch functional derived from the given
boundary value problem. Several versions of the convexification method have
been developed since it was first introduced in \cite%
{KlibanovIoussoupova:SMA1995}. We cite here \cite%
{Klibanov:ip2015,Klibanov:sjma1997,KlibanovNik:ra2017,KlibanovLiZhang:ip2019,KhoaKlibanovLoc:SIAMImaging2020,KlibanovLiZhang:SIAM2019}
for some important works in this area and their applications to solve a
variety kinds of inverse problems.
A comprehensive study of the convexification method is presented in the recent published book \cite{KlibanovLiBook}.
The crucial mathematical results that
guarantee the above mentioned properties of the convexification, are the
Carleman estimates. The original idea of applying Carleman estimates to
coefficient inverse problems was first published in \cite%
{BukhgeimKlibanov:smd1981} back in 1981 to prove uniqueness theorems for a
wide class of coefficient inverse problems. Some follow up publications can
be found in, e.g. \cite{KlibanovTimonov:u2004,Klibanov:ip1992,Isakov:bookSpringer2017,LocNguyen:ip2019,LeNguyen:2020,TriggianiandYao:amo2002}%
. Surveys on the method in \cite{BukhgeimKlibanov:smd1981} can be found in 
\cite{Klibanov:jiipp2013,Yamamoto:ip2009}, see also \cite[Chapter 1]%
{BeilinaKlibanovBook} and \cite{KlibanovLiBook}. It was discovered later in \cite%
{KlibanovIoussoupova:SMA1995}, that the idea of \cite%
{BukhgeimKlibanov:smd1981} can be successfully modified to develop globally
convergent numerical methods for coefficient inverse problems using the
convexification.

For the convenience of the reader, we will recall in this paper the
convexification method , to solve ill-posed Cauchy problems for quasi-linear
PDEs with both Dirichlet and Neumann boundary data \cite{Klibanov:ip2015}. Then, we will prove that if the noise in the boundary
data tends to zero, then the convexification method combined with the
gradient descent method delivers a close approximation to the solution of
that Cauchy problem if starting from an arbitrary point of a selected
bounded set. The rate of convergence is Lipschitz. We next apply the above
results to solve a highly nonlinear and severely ill-posed coefficient
inverse problem, described below. At a point far away from the region of
interest, we send out an incident electric wave. The incident wave
propagates in the 3D space and scatters when hitting the targets. We measure
the back scattering wave on a surface. The aim of the inverse problem is to
reconstruct the spatially distributed dielectric constants from this
measurement. This coefficient inverse problem is the so-called inverse back
scattering problem. It has many real-world applications, including the
detection and identification of explosives, nondestructive testing and
material characterization, see \cite%
{Khoaelal:IPSE2021,VoKlibanovNguyen:IP2020,Schubert:sp2006,Weatherall:2015}.
We also refer the reader to \cite{ColtonKress:2013} for the applications of
this coefficient inverse problem in sonar imaging, geographical exploration,
medical imaging, near-field optical microscopy, nano-optics.

The widely-used method to solve nonlinear coefficient inverse problem is the
least squares optimization. This approach requires a good initial guess of
the true solution. 
Unlike this, we assume that the target to be detected is completely
unknown. This means that a good initial guess of its dielectric constant is
unavailable. Our numerical procedure is as follows. We first eliminate the
unknown dielectric constant from the governing Helmholtz equation. The
obtained equation is not a standard PDE. Then, we approximate the solution
of the latter PDE via a truncated Fourier series with respect to a special
orthonormal basis. Then we obtain an ill-posed Cauchy problem for a coupled
system of elliptic of PDEs wih respect to corresponding spatially dependent
Fourier coefficients. That special orthonormal basis was originally
introduced in \cite{Klibanov:jiip2017}. Solving this system by the
convexification method and the gradient descent method, we obtain the
solution to the above non standard PDE above. Then the solution to the
originating coefficient inverse problem follows. We refer the reader to \cite%
{VoKlibanovNguyen:IP2020,KhoaKlibanovLoc:SIAMImaging2020,TruongNguyenKlibanov:IPSE2021}
and the references therein for some related versions of this method.

The paper is organized as follows. In Section \ref{convex ana}, we prove the
convergence of the gradient descent method to the minimizer of a strictly
convex functional. In Section \ref{sec statement 1}, we present the above
mentioned ill-posed Cauchy problem for a quasilinear elliptic PDE with both
Dirichlet and Neumann boundary conditions. Also, in this section we present
the corresponding functional with the Carleman Weight Function in it. In
section \ref{sec convex}, we recall the convexification method. In this
section, we also prove the Lipschitz-like convergence of the minimizers due
to the convexification method to the true solution as the noise tends to
zero. 
In section \ref{sec CIP}, we introduce our coefficient inverse problem. In
section \ref{sec 7}, we derive an approximate model to solve this
coefficient inverse problem. We present some numerical examples in section %
\ref{sec num}. Section \ref{sec 9} is for concluding remarks.

\section{The gradient descent method to minimize a convex functional}

\label{convex ana}

Let $X$ be a Hilbert space and let $J:X\rightarrow \mathbb{R}$ be a
functional. 
Assume that $J$ is Fr\'echet differentiable. Its derivative at the point $%
v\in X$ is denoted by $DJ(v):X\rightarrow \mathbb{R}$. By the Riesz
representation theorem, for each $v\in X,$ we can identify $DJ(v)$ with an
element of $X$, named as $J^{\prime }(v),$ in the following sense 
\begin{equation*}
DJ(v)(h)=\langle J^{\prime }(v),h\rangle _{X}\quad \mbox{for all }h\in X.
\end{equation*}%
Let $M>0$ be an arbitrary number. Consider the ball $
B\left( M\right) ,$%
\begin{equation*}
B\left( M\right) =\left\{ v\in X:\left\Vert v\right\Vert _{X}<M\right\} .
\end{equation*}%
 We assume that the Fr\'echet derivative $J^{\prime }$ is Lipschitz continuous in $\overline{B\left(
M\right)} ,$ i.e.%
\begin{equation}
\left\Vert J^{\prime }(v_{2})-J^{\prime }(v_{1})\right\Vert _{X}\leq
L\left\Vert v_{2}-v_{1}\right\Vert _{X},\text{ }\forall v_{1},v_{2}\in
\overline{B\left(
M\right)} ,  \label{1}
\end{equation}%
where $L$ is a certain number. 
Assume that $J$ is strictly convex on $\overline{B\left(
M\right)} $. This means that there exists a constant $\Lambda >0$
such that 
\begin{equation}
\langle J^{\prime }(v_{1})-J^{\prime }(v_{2}),v_{1}-v_{2}\rangle _{X}\geq
\Lambda \Vert v_{1}-v_{2}\Vert _{X}^{2}\quad \mbox{for all }v_{1},v_{2}\in
\overline{B\left(
M\right)}.  \label{convexity1}
\end{equation}

Theorem \ref{thm 2.11} follows from a combination of Lemma 2.1 and Theorem 2.1 of 
\cite{KlibanovNik:ra2017}.

\begin{theorem}
Assume that the functional $J:X\rightarrow \R$
 is Fr\'echet differentiable on $X $
and its Fr\'echet derivative is Lipschitz continuous on $\overline{B(M)}$ as in 
\eqref{1}. Also, assume that $J\left( v\right) $ is  strictly convex
in $\overline{B(M)}$; i.e., inequality \eqref{convexity1}
is true. Then there exists unique minimizer $v_{\mathrm{min}}\in \overline{%
B\left( M\right) }$ of the functional $J\left( v\right) $ on
the set $\overline{B\left( M\right) },$%
\begin{equation*}
\min_{\overline{B\left( M\right) }}J\left( v\right) =J\left( v_{\mathrm{min}%
}\right) .
\end{equation*}%
Furthermore, the following inequality holds
\begin{equation*}
\langle J^{\prime }\left( v_{\mathrm{min}}\right) ,v_{\mathrm{min}}-y\rangle_X
\leq 0, \quad \mbox{for all } y\in \overline{B\left( M\right) }.
\end{equation*}
\label{thm 2.11}
\end{theorem}

The fact that the minimizer guaranteed by Theorem \ref{thm 2.11} can be located
on the boundary of the ball $\overline{B(M)}$
prevents us from the proof of the global convergence of the gradient descent
method. Hence, we assume in the next theorem that the minimizer is an
interior point of $B(M).$

\begin{theorem}
Assume that the functional $J:X\rightarrow \mathbb{R}$ satisfies conditions
of Theorem \ref{thm 2.11}. Let $v_{\mathrm{min}}$ be its minimizer on the set $\overline{%
B\left( M\right) },$ the uniqueness and existence of which is guaranteed by
Theorem \ref{thm 2.11}. Suppose that $v_{\mathrm{min}}$ belongs to the interior of $%
B\left( M\right) .$ Fix $v^{(0)}\in B$. Assume that the ball centered at $v_{%
\mathrm{min}}$ with the radius $\Vert v^{(0)}-v_{\mathrm{min}}\Vert _{X}$ is
contained in $B\left( M\right);$ i.e.,
\begin{equation}
B_{0}=B(v_{\mathrm{min}},\Vert v^{(0)}-v_{\mathrm{min}}\Vert _{X})\subset
B\left( M\right).
\label{interior}
\end{equation}%
Denote $\eta _{0}=\min \left( 2\Lambda /L^{2},1\right)$ and fix $\eta \in (0, \eta_0).$ 
For each $m \geq 0,$ define  
\begin{equation}
v^{(m+1)}=v^{(m)}-\eta J^{\prime }\big( v^{(m)}\big) ,\quad m\geq 1
\label{minimizing sequence}
\end{equation}%
Then, there exists a number $q\in \left( 0,1\right) 
$ such that 
\begin{equation}
v^{(m)} \in B(M) \quad
\mbox{and }
\big\Vert v^{m}-v_{\min }\big\Vert _{X} \leq
 q^{m-1}\big\Vert
v^{(0)}-v_{\min }\big\Vert _{X}, \quad m \geq 1.  \label{2}
\end{equation}
As a result, the sequence $v^{(m)}$ converges to $v_{\rm min}$ as $m$ tends to $\infty.$
\label{gdm}
\end{theorem}

The sequence $\{v^{(m)}\}_{m\geq 1}$ defined in \eqref{minimizing sequence}
is generated by the well-known gradient descent method. Although the
gradient descent method is widely used in the scientific community,
its convergence for a nonconvex functional can be proven only if the
starting point of iterations is sufficiently close to the minimizer. Unlike
this, Theorem \ref{gdm} provides an affirmative answer about the convergence
when the starting point is not necessary located in a sufficiently small
neighborhood of the minimizer. Theorem \ref{gdm} justifies recent numerical
results of our research group where we have used the gradient descent method
 to  minimize globally strictly convex cost functionals arising in
convexification even though our past theory said that a more complicated
gradient projection method should be used, see e.g., \cite%
{KlibanovNik:ra2017,VoKlibanovNguyen:IP2020,Khoaelal:IPSE2021,KhoaKlibanovLoc:SIAMImaging2020,KlibanovLiBook,KlibanovLiZhang:SIAM2019}. 

\textit{Proof of Theorem \ref{gdm}.} Let $L$ and $\Lambda $ be the constants
in \eqref{1} and \eqref{convexity1} respectively and let $q = 1+\eta ^{2}L^{2}-2\eta\Lambda$. 
Since $\eta \in (0, \eta_0),$ $q \in (0, 1).$ We prove \eqref{2} by induction.
 Assume, by
induction, that \eqref{2} is true for some $m\geq 1$. 
Due to assumption \eqref{interior}, 
$v_{\min }\in B\left( M\right).$ Thus, $J^{\prime }(v_{\min })=0.$ Hence, $%
v_{\min }=v_{\min }-J^{\prime }\left( v_{\min }\right) .$ By %
\eqref{minimizing sequence}, we have 
\begin{align*}
\Vert v^{(m+1)}-v_{\min }\Vert _{X}^{2}& =\Vert v^{(m)}-v_{\min }-\eta
(J^{\prime }\left( v^{(m)}\right) -J^{\prime }(v_{\min }))\Vert _{X}^{2} \\
& =\Vert v^{(m)}-v_{\min }\Vert _{X}^{2}+\eta ^{2}\Vert J^{\prime }\left(
v^{(m)}\right) -J^{\prime }(v_{\min })\Vert _{X}^{2} \\
& \hspace{2.5cm}-2\eta \langle J^{\prime (m)})-J^{\prime }(v_{\min
}),v^{(m)}-v_{\min }\rangle _{X}.
\end{align*}%
Using this, together with \eqref{1}, \eqref{convexity1} and the induction assumption for \eqref{2}, we obtain 
\begin{equation}
\Vert v^{(m+1)}-v_{\min }\Vert _{X}^{2}\leq q\Vert v^{(m)}-v_{\min }\Vert _{X}^{2}
\leq q^{m}\Vert v^{(0)}-v_{\min }\Vert _{X}^{2}.  \label{4}
\end{equation}
The last inequality in \eqref{4} is deduced from the induction hypothesis.
It follows from \eqref{4} that $ v^{(m+1)} \in B_0 
\subset B$. The assertion \eqref{2} is proved.
 \QEDB%\end{proof}

\begin{Remark}
\begin{enumerate}
\item The hypothesis that the starting point of iterations $v^{(0)}$ is such
that the ball centered at $v_{\mathrm{min}}$ with the radius $\Vert v_{\min
}-v_{0}\Vert _{X}$ is contained in $B\left( M\right) $ does not weaken
Theorem \ref{gdm}. In fact, if this hypothesis is not satisfied, we can
replace $B\left( M\right) $ by a larger ball $B\left( M^{\prime }\right) $ 
where $M^{\prime }>M.$ Note that in the convexification method, $B\left(
M\right) $ is the ball with an arbitrary chosen radius. 

\item The assumption that $v_{\mathrm{min}}$ is inside $B\left( M\right) $
is the main reason that helps us to replace the gradient projection method
in \cite%
{Khoaelal:IPSE2021,KhoaKlibanovLoc:SIAMImaging2020,VoKlibanovNguyen:IP2020,KlibanovNik:ra2017}
with the gradient descent method in Theorem \ref{gdm}. Without this
assumption, elements of the sequence produced by the gradient descent method
might be outside of $B\left( M\right) $, thus making this sequence diverge. 

\item We refer the reader to \cite[Theorem 6]{Klibanov:2ndSAR2021}, in which the authors proved a less general case of Theorem \ref{gdm}.

\end{enumerate}
\end{Remark}

\section{A boundary value problem for quasi-linear PDEs}

\label{sec statement 1} Let $n\geq 2$ be the spatial dimension. Let $\Omega $
be an open and bounded domain in $\mathbb{R}^{n}$ and $\Gamma $ be a part of 
$\partial \Omega .$ Let $G:\Omega \times \mathbb{R}\times \mathbb{R}%
^{n}\rightarrow \mathbb{R}$ be a real value function in the class $C^{2}(%
\overline{\Omega }\times \mathbb{R}\times \mathbb{R}^{n},\mathbb{R})$.
Consider the following boundary value problem with both Dirichlet and
Neumann boundary conditions 
\begin{equation}
\left\{ 
\begin{array}{rcll}
\Delta v(\mathbf{x}) & = & G(\mathbf{x},v(\mathbf{x}),\nabla v(\mathbf{x}))
& \mathbf{x}\in \Omega , \\ 
\partial _{\nu }v(\mathbf{x}) & = & g_{0}(\mathbf{x}) & \mathbf{x}\in \Gamma
, \\ 
v(\mathbf{x}) & = & g_{1}(\mathbf{x}) & \mathbf{x}\in \partial \Omega 
\end{array}%
\right.   \label{nonlinear eqn}
\end{equation}%
where $g_{0}$ and $g_{1}$ are two functions in the class $H^{p}(\Omega )$
where $p$ is a positive integer with $p>\ceil{n/2}+2.$ In fact, we
can say that \eqref{nonlinear eqn} is the Cauchy problem for a quasilinear
elliptic equation with the additional Dirichlet boundary data at $\partial
\Omega \diagdown \Gamma .$ Here, $\ceil{n/2}$ is the smallest integer that
is greater than $n/2$. This regularity condition guarantees the embedding $%
H^{p}(\Omega )\hookrightarrow C^{2}(\overline{\Omega })$. This embedding
will be used for the regularization purpose. In practice, the functions $%
g_{0}$ and $g_{1}$ represent the flux and the value information of $v$ on $%
\partial \Omega $ and $\Gamma $ respectively. We first recall the
convexification method to compute an approximation of the solution, if
exists, to \eqref{nonlinear eqn}. Suppose that 
\begin{equation}
H=\{\phi \in H^{p}(\Omega ):\partial _{\nu }\phi (\mathbf{x})=g_{0}(\mathbf{x%
}) 
\mbox{ for all }\mathbf{x}\in \Gamma \mbox{ and }\phi (\mathbf{x})=g_{1}(%
\mathbf{x})\mbox{ for all }\mathbf{x}\in \partial \Omega \}  \label{2,2}
\end{equation}%
is nonempty. Let $v_{0}$ be a function in $H$. Define 
\begin{equation}
u(\mathbf{x})=v(\mathbf{x})-v_{0}(\mathbf{x})\quad \mbox{for all }\mathbf{x}%
\in \Omega .  \label{change}
\end{equation}%
Then, solving \eqref{nonlinear eqn} is equivalent to solving 
\begin{equation}
\left\{ 
\begin{array}{rcll}
\Delta u(\mathbf{x}) & = & F(\mathbf{x},u(\mathbf{x}),\nabla u(\mathbf{x}))
& \mathbf{x}\in \Omega , \\ 
\partial _{\nu }u(\mathbf{x}) & = & 0 & \mathbf{x}\in \Gamma , \\ 
u(\mathbf{x}) & = & 0 & \mathbf{x}\in \partial \Omega 
\end{array}%
\right.   \label{nonlinear v}
\end{equation}%
where 
\begin{equation}
F(\mathbf{x},s,\xi )=\Delta v_{0}(\mathbf{x})+G(\mathbf{x},s+v_{0}(\mathbf{x}%
),\xi +\nabla v_{0}(\mathbf{x}))  \label{change F}
\end{equation}%
for all $\mathbf{x}\in \Omega ,s\in \mathbb{R},\xi \in \mathbb{R}^{n}$. Let 
\begin{equation}
H_{0}=\big\{\phi \in H^{p}(\Omega ):\partial _{\nu }\phi (\mathbf{x})=0 
\mbox{ for all }\mathbf{x}\in \Gamma \mbox{ and }\phi (\mathbf{x})=0\mbox{
for all }\mathbf{x}\in \partial \Omega \big\}.  \label{H0}
\end{equation}%
It is obvious that $H_{0}$ is a closed subspace of $H^{p}(\Omega )$. We
consider $H_{0}$ a Hilbert space endowed with the usual norm of $%
H^{p}(\Omega ).$ 
A widely-used approach to solve \eqref{nonlinear v} is to minimize the
following least squares functional 
\begin{equation}
\int_{\Omega }\big|\Delta u(\mathbf{x})-F\big(\mathbf{x},u(\mathbf{x}%
),\nabla u(\mathbf{x})\big)\big|^{2}d\mathbf{x}+\mbox{a regularization term}.
\label{nonconvex}
\end{equation}%
for $u\in H_{0}.$ Due to the nonlinearity of $F$, and hence $L$, the
functional in \eqref{nonconvex} is nonconvex. It might have multiple local
minima and ravines, making the direct optimization approach unpractical.
Motivated by this fact, we \textquotedblleft convexify" this functional
using the idea in \cite{Klibanov:jiip2017}. Let $\Psi :\Omega \rightarrow 
\mathbb{R}$ be a $C^{\infty }$ function with $\nabla \Psi \not=0$ for all $%
\mathbf{x}\in \overline{\Omega }.$ Introduce the Carleman weight function 
\begin{equation}
\mu _{\lambda }(\mathbf{x})=\exp (2\lambda \Psi (\mathbf{x}))\quad \mbox{for
all }\mathbf{x}\in \Omega   \label{Carleman weight}
\end{equation}%
where $\lambda >1$. A Carleman estimate is an inequality of the form below.

\begin{Assumption}[Carleman estimate]
\textit{There exists $\lambda_0 > 1$ depending only on $\Omega$ and $n$ such
that for all $\phi \in H^2(\Omega)$ with $\phi|_{\partial \Omega} = 0$ and $%
\partial_\nu \phi|_{\Gamma} = 0$, we have  
\begin{equation}
\int_{\Omega} \mu_\lambda(\mathbf{x}) |\Delta \phi|^2 d\mathbf{x} \geq \frac{%
C}{\lambda} \sum_{i, j = 1}^n \int_{\Omega}\mu_\lambda(\mathbf{x}%
)|\partial_{x_i x_j}\phi(\mathbf{x})|^2d\mathbf{x} 
+ C \lambda \int_{\Omega} \mu_\lambda(\mathbf{x}) [|\nabla \phi|^2 +
\lambda^2 |\phi|^2] d\mathbf{x}  \label{Car est}
\end{equation}
for all $\lambda \geq \lambda_0$ for some positive constant $C$ depending
only on $\Omega$. \label{assumption Carleman} }
\end{Assumption}

Assumption \ref{assumption Carleman} holds true for some functions $\Psi$
and $\mu_\lambda$. For example, Klibanov and his collaborators have
established a Carleman estimate in \cite[Theorem 4.1]{KlibanovLiZhang:ip2019}%
, in which $\Omega = (-R, R)^3 \subset \mathbb{R}^3$ for some $R > 0$ and $%
\Psi(x, y, z) = (z - r)^2$ where $r$ is any number that is greater than $R$.
On the other hand, following the arguments in \cite[Theorem 4.1]%
{KlibanovLiZhang:ip2019}, one can prove a Carleman estimate with $\Psi(x, y,
z) = (z + r)^2$, see \cite[Theorem 3.1]{KhoaKlibanovLoc:SIAMImaging2020}.
This estimate plays an important role in developing a numerical method to
solve the back scattering inverse problem with moving point source in \cite%
{KhoaKlibanovLoc:SIAMImaging2020}. On the other hand, the reader can find
another Carleman estimate in \cite[Theorem 3.1]{LeNguyen:2020} when the second derivatives of the
test function $\phi$ are absence in the right hand side of \eqref{Car est}.
We also cite to \cite{BeilinaKlibanovBook,KlibanovLiBook} for some important versions of
Carleman estimate for other kinds of partial differential operators.
Especially, we draw the reader's attention to \cite{BukhgeimKlibanov:smd1981}
for the original idea of using Carleman estimate to prove the uniqueness of
a variety kinds of inverse problems.

Without lost of the generality, we assume that the true solution to \eqref{nonlinear v} has a finite $H^p(\Omega)$-norm, which is bounded from above by a known number $M$.
We  seek this solution in the set 
\begin{equation}
B(M)=\big\{\phi \in H_{0}:\Vert \phi \Vert _{H^{p}(\Omega )}<M\big\}.
\label{bm}
\end{equation}
More precisely, in order to find a numerical solution to \eqref{nonlinear v}, we solve the
problem below.

\begin{problem}
\textit{\ Fix a regularization parameter $\epsilon > 0$. Minimize the
following functional  
\begin{equation}
J_{\epsilon, \lambda}(u) = \int_\Omega \mu_\lambda(\mathbf{x}) \big|\Delta u(%
\mathbf{x}) - F\big(\mathbf{x}, u(\mathbf{x}), \nabla u(\mathbf{x})\big)\big|%
^2 d\mathbf{x} + \epsilon \|u\|^2_{H^p(\Omega)}  \label{objective function}
\end{equation}
for all $u \in B(M)$.  \label{pro optimization} }
\label{p1}
\end{problem}

%\begin{Remark}
In the next section, we will recall the convexification principle to solve
Problem \ref{pro optimization}.  The main content of the convexification
principle is that if the Carleman estimate \eqref{Car est} holds true then
for any arbitrarily large number $M$, there exists $\lambda_1 > \lambda_0$
such that for all $\lambda > \lambda_1$ and $\epsilon > 0,$ $J_{\epsilon,
\lambda}$ is strictly convex in $\overline{B(M)}$ where $B(M)$ is the ball
in $H_0$ with center $0$ and radius $M$, see \eqref{bm}.

\section{The convexification method and the convergence of the minimizer to
the true solution as the noise tends to zero}

\label{sec convex}

% For $M>0$, define 
% \begin{equation}
% B(M)=\big\{\phi \in H_{0}:\Vert \phi \Vert _{H^{p}(\Omega )}<M\big\}.
% \label{bm}
% \end{equation}%
In this section, we recall a theorem (Theorem \ref{thm 2.1}) that guarantees
the convexity of the objective functional $J_{\epsilon ,\lambda }$ in $%
\overline{B(M)}$. We write $F=F(\mathbf{x},s,\xi )$ and the partial
derivatives of $F(\mathbf{x},s,\xi )$ with respect to its variables are
written as $\nabla _{\mathbf{x}}F(\mathbf{x},s,\xi ),$ $\partial _{s}F(%
\mathbf{x},s,\xi )$ and $\nabla _{\xi }F(\mathbf{x},s,\xi )$. The following
theorem, Theorem \ref{thm 2.1}, guarantees that $J_{\epsilon ,\lambda }$ has
Lipschitz continuous Fr\'echet derivative and, more importantly, that $%
J_{\epsilon ,\lambda }$ is strictly convex if the Carleman weight function $%
\mu _{\lambda }$ is such that Assumption \ref{assumption Carleman} holds
true.

\begin{theorem}[Convexification]
1. Let $M$ be an arbitrary positive number and define the ball $B(M)$ as in %
\eqref{bm}.  Then, for all $\epsilon > 0$ and $\lambda > 1$, $J_{\epsilon,
\lambda}: \overline{B(M)} \subset H_0 \to \mathbb{R}$ is Fr\'echet
differentiable. The derivative of $J_{\epsilon, \lambda}$ is given by 
\begin{equation}
DJ_{\epsilon, \lambda}(u)h = 2\int_{\Omega} \mu_\lambda(\mathbf{x}) \big(%
\Delta u(\mathbf{x}) - F(\mathbf{x}, u(\mathbf{x}), \nabla u(\mathbf{x}))%
\big) 
\big(\Delta h(\mathbf{x}) - DF(u)h(\mathbf{x})\big) d\mathbf{x} +2 \epsilon
\langle u, h \rangle_{H^p(\Omega)}  \label{DJ}
\end{equation}
for all $u \in {B(M)}$ and $h \in H_0$ where 
\begin{equation*}
DF(u)h(\mathbf{x}) = \partial_s F(\mathbf{x}, u(\mathbf{x}), \nabla u(%
\mathbf{x})) h(\mathbf{x}) + \nabla_{\xi}F(\mathbf{x}, u(\mathbf{x}), \nabla
u(\mathbf{x}))\cdot \nabla h(\mathbf{x}) 
\end{equation*}
for all $\mathbf{x} \in \Omega.$ Moreover,  the Fr\'echet derivative $%
DJ_{\epsilon, \lambda}$  is Lipschitz continuous in $\overline{B(M)}$.  That
means,  there exists a constant $L = L(\Omega, M, F)$, depending only on the
listed parameters, such that  
\begin{equation}
\|DJ_{\epsilon, \lambda}(u_2) - DJ_{\epsilon, \lambda}(u_1)\|_{\mathcal{L}%
(H_0)} \leq L\|u_2 - u_1\|_{H^p(\Omega)}  
%\quad \mbox{for all } u_1, u_2 \in \overline{B(M)}
\label{DJLip}
\end{equation}
for all $u_1, u_2 \in \overline{B(M)}$,
where $\mathcal{L}(H_0)$ is the set of all bounded linear maps sending
functions in $H_0$ into $\mathbb{R}.$

2. Assume further that the Carleman estimate \eqref{Car est} holds true. 
Then, there exist $\lambda_1 = \lambda(M, \Omega, F)> \lambda_0$ and $C =
C(M, \Omega, F) > 0$, both of which depend only on the listed parameters,
such that for all $\epsilon > 0,$ $\lambda > \lambda_1$,  $u_1$ and $u_2$ in 
$\overline{B(M)}$, we have  
\begin{equation}
J_{\epsilon, \lambda}(u_2) - J_{\epsilon, \lambda}(u_1) - DJ_{\epsilon,
\lambda}(u_1)(u_2 - u_1) \\
\geq C\|u_2 - u_1\|_{H^2(\Omega)}^2 + \epsilon \|u_2 - u_1\|^2_{H^p(\Omega)}.
\label{convexity}
\end{equation}
As a result,  
\begin{equation}
\Big(DJ_{\epsilon, \lambda}(u_2) - DJ_{\epsilon, \lambda}(u_1) \Big)(u_2 -
u_1) \\
\geq 2C\|u_2 - u_1\|^2_{H^2(\Omega)} + 2\epsilon\|u_2 - u_1\|^2_{H^p(\Omega)}
\label{3,2}
\end{equation}
for all $u_1$ and $u_2$ in $\overline{B(M)}$.

3. $J_{\epsilon, \lambda}$ has a unique minimizer in $\overline{B(M)}.$  %
\label{thm 2.1}
\end{theorem}

We do not present the proof of Theorem \ref{thm 2.1} here. The reason is below.
One can prove the first part of this theorem with straight forward
computations. The proof is similar to that of \cite[Theorem 3.1]%
{KlibanovNik:ra2017}. The second part of this theorem is a generalization of 
\cite[Theorem 3.2]{KlibanovNik:ra2017} in the sense that the
\textquotedblleft convexification" inequalities \eqref{convexity} and %
\eqref{3,2} are tighter than the ones in \cite[Theorem 3.2]%
{KlibanovNik:ra2017}. In fact, in those inequalities, we replace the $H^{1}$
norm in \cite[Theorem 3.2]{KlibanovNik:ra2017} by the $H^{2}$ norm in the
right hand side of \eqref{convexity} and \eqref{3,2}. This is because the
right hand side of the Carleman estimate, see \eqref{Car est}, contains the
second derivatives. The existence of the unique minimizer of $J_{\epsilon
,\lambda }$ in part 3 of Theorem \ref{thm 2.1} can be proved using the same
technique in \cite[Theorem 5.3]{KhoaKlibanovLoc:SIAMImaging2020}, see also Lemma 2.1 and Theorem 2.1 of 
\cite{KlibanovNik:ra2017} and Theorem \ref{thm 2.11}. 
On the
other hand, we refer the reader to \cite[Section 2]{KlibanovNik:ra2017} for
some important facts in convex analysis that are related to the
convexification in Theorem \ref{thm 2.1}.

By using the gradient descent method, we can compute the minimizer of $%
J_{\epsilon, \lambda}$ in $B(M)$, see Theorem \ref{gdm}. We are now in the
position of solving problem \eqref{nonlinear eqn} with noisy boundary data $%
g_0$ and $g_1$ given. The corresponding noiseless data are denoted by $g_0^*$
and $g_1^*$ respectively. Let $\delta > 0$ be the noise level and assume
that there exists an ``error function" $\mathcal{E}$ such that 
\begin{equation}
\left\{ 
\begin{array}{ll}
\|\mathcal{E}\|_{H^p(\Omega)} < \delta, &  \\ 
g_0 = g_0^* + \partial_{\nu }\mathcal{E} & \mbox{on } \partial \Omega, \\ 
g_1 = g_1^* + \mathcal{E} & \mbox{in } \Gamma.%
\end{array}
\right.  \label{19}
\end{equation}
Recall in Section \ref{sec statement 1}, we assume that there is a function $%
v_0$ satisfying $\partial_{\nu} g = g_0$ on $\partial \Omega$ and $v_0 = g_1$
on $\Gamma$. Let 
\begin{equation}
v_{\epsilon, \delta}(\mathbf{x}) = u_{\min}(\mathbf{x}) + v_0(\mathbf{x})
\quad \mbox{for all } \mathbf{x} \in \Omega  \label{4.8}
\end{equation}
where $u_{\min}$ is the minimizer of $J_{\epsilon, \lambda}$ obtained in
Theorem \ref{thm 2.1}. The function $v_{\epsilon, \delta}$ is named as the
regularized solution to \eqref{nonlinear eqn}. Let $v^*$ be the solution to %
\eqref{nonlinear eqn} with $g_0$ and $g_1$ replaced by the corresponding
noiseless data $g_0^*$ and $g_1^*$ respectively.  The following theorem
confirms that the minimizer of $J_{\epsilon, \lambda}$ can be used to
approximate the solution to \eqref{nonlinear eqn} via \eqref{4.8}.  It is a
generalization of Theorem 4.5 in \cite{KhoaKlibanovLoc:SIAMImaging2020} and
Theorem 5.4 in \cite{KhoaKlibanovLoc:SIAMImaging2020}. In fact, in those
theorems, the function $F$ has some specific form and does not depend on the
first and the second variables $\mathbf{x}$ and $u(\mathbf{x})$. 
% Its proof follows closely to the proof of Theorem 4.5 in \cite{KhoaKlibanovLoc:SIAMImaging2020} and Theorem 5.4 in \cite{KhoaKlibanovLoc:SIAMImaging2020}.

\begin{theorem}
Assume that problem \eqref{nonlinear eqn} with $g_0$ and $g_1$ replaced by $%
g_0^*$ and $g_1^*$ respectively has a solution $v^*$. Recall $v_0$ the
function we used to change the variable in \eqref{change}.  Without loss of
the generality, assume that  
\begin{equation}
\max\big\{\|v^*\|_{H^p(\Omega)}, \|v_0\|_{H^p(\Omega)}\big\} < \frac{M}{2} -
\delta.  \label{v*M}
\end{equation}
Let $v_{\epsilon, \delta} = u_{\min} + v_0$ where $u_{\min}$ is the
minimizer of the strictly convex functional $J_{\epsilon, \lambda}$.  Then  
\begin{equation}
\|v_{\epsilon, \delta} - v^*\|_{H^2(\Omega)} \leq C(\sqrt{\epsilon}\|v^* -
v_0\|_{H^p(\Omega)} + \delta)  \label{stable V}
\end{equation}
for some constant $C$ depending on $\Omega, M$ and $F$.  \label{thm 5.1}
\end{theorem}

%\begin{proof}
\noindent\textit{Proof. }  For $\mathbf{x} \in \Omega$, define  
\begin{equation}
u^*(\mathbf{x}) = v^*(\mathbf{x}) - v_0(\mathbf{x}).  \label{verror}
\end{equation}
It is obvious that  
\begin{equation*}
\left\{  
\begin{array}{ll}
\partial_{\nu} u^*(\mathbf{x}) = g_0^*(\mathbf{x}) - g_0(\mathbf{x}) = -
\partial_{\nu}\mathcal{E}(\mathbf{x}) & \mbox{for all }\mathbf{x} \in
\partial \Omega, \\ 
u^*(\mathbf{x}) = g_1^*(\mathbf{x}) - g_1(\mathbf{x}) = - \mathcal{E}(%
\mathbf{x}) & \mbox{for all }\mathbf{x} \in \Gamma%
\end{array}
\right.  
\end{equation*}
Thus, $u^* + \mathcal{E }\in H_0$.  Using the triangle inequality, \eqref{19}
and \eqref{v*M}, we have $u^* + \mathcal{E }\in B(M).$  Using %
\eqref{convexity} with $u_1$ and $u_2$ replaced by $u_{\mathrm{min}}$ and $%
u^* + \mathcal{E}$ respectively, we have  
\begin{multline}
J_{\epsilon, \lambda}(u^* + \mathcal{E}) - J_{\epsilon, \lambda}(u_{\mathrm{%
min}}) - DJ_{\epsilon, \lambda}(u_{\mathrm{min}})(u^* + \mathcal{E }-
u_{\min}) \\
\geq C\|u^* + \mathcal{E }- u_{\min}\|_{H^2(\Omega)}^2 + \epsilon \|u^* + 
\mathcal{E }- u_{\min}\|^2_{H^p(\Omega)}.  \label{22}
\end{multline}
Since $u_{\min}$ is the minimizer of $J_{\epsilon, \lambda}$ in $B(M)$, $%
DJ_{\epsilon, \lambda}(u_{\mathrm{min}}) = 0$. This, together with \eqref{22}
and the fact that $- J_{\epsilon, \lambda}(u_{\mathrm{min}}) \leq 0$,
implies 
\begin{equation}
J_{\epsilon, \lambda}(u^* + \mathcal{E}) \geq C\|u^* + \mathcal{E }-
u_{\min}\|_{H^2(\Omega)}^2 + \epsilon \|u^* + \mathcal{E }-
u_{\min}\|^2_{H^p(\Omega)}.  \label{24}
\end{equation}
Using the inequality $(a + b)^2 \leq 2a^2 + 2b^2,$ we next estimate 
\begin{align}
J_{\epsilon, \lambda}&(u^* + \mathcal{E})\nonumber 
= \int_{\Omega} \mu_\lambda(%
\mathbf{x}) |\Delta (u^* + \mathcal{E}) - F(\mathbf{x}, u^* + \mathcal{E},
\nabla u^* + \nabla \mathcal{E}) |^2d\mathbf{x} 
+\epsilon \|u^* + \mathcal{E}%
\|_{H^p(\Omega)}^2  \notag \\
&\leq 2\int_{\Omega} \mu_\lambda(\mathbf{x}) |\Delta u^* - F(\mathbf{x},
u^*, \nabla u^*) |^2d\mathbf{x}  \notag \\
&\quad + 2\int_{\Omega} \mu_\lambda(\mathbf{x}) |\Delta \mathcal{E }+ F(%
\mathbf{x}, u^*, \nabla u^*) - F(\mathbf{x}, u^* + \mathcal{E}, \nabla u^* +
\nabla \mathcal{E}) |^2d\mathbf{x} 
%\notag
%\\
%&\hspace{3cm}
+\epsilon \|u^* + \mathcal{E}%
\|_{H^p(\Omega)}^2.    \label{28}
\end{align}
Since $v^*$ is the true solution to \eqref{nonlinear eqn}, by %
\eqref{nonlinear v}, \eqref{change F} and \eqref{verror}, we have 
\begin{equation*}
\Delta u^* - F(\mathbf{x}, u^*, \nabla u^*) = 0 \quad \mbox{for all } 
\mathbf{x} \in \Omega.
\end{equation*}
Using \eqref{19} and \eqref{28} and the fact that $F$ is in $C^1$ and
hence Lipschitz, we have 
\begin{equation}
J_{\epsilon, \lambda}(u^* + \mathcal{E}) \leq C \delta^2 + \epsilon \|u^* + 
\mathcal{E}\|_{H^p(\Omega)}^2.  \label{25}
\end{equation}
Combining \eqref{24} and \eqref{25} and using the inequality $(a + b)^2 \leq
2(a^2 + b^2)$, we have 
\begin{equation*}
\|u^* - u_{\min}\|_{H^2(\Omega)}^2 \leq C (\delta^2 + \epsilon
\|u^*\|^2_{H^p(\Omega)}). 
\end{equation*}
Estimate \eqref{stable V} is proved.  \QEDB

The procedure to compute $v^*$ is described in Algorithm \ref{alg 1}.

\begin{algorithm}
	\caption{\label{alg 1} A numerical method to solve \eqref{nonlinear eqn}}
	\begin{algorithmic}[1]
	\State \label{choice} Choose $M$ large enough and choose a threshold error $\varepsilon > 0$.  
	\State  \label{choose U0} Set $m = 0$ and take a function $u_0$ in $B(M)$. 
	\State \label{step backward} Compute $u_{m + 1}$ using the gradient descent method, see \eqref{minimizing sequence} for some $0 < \eta \ll 1$. 
	\State \label{step 4} If $\|u_{m + 1} - u_{m}\|_{H^2(\Omega)} < \varepsilon$, go to Step \ref{step forward}. Otherwise, set $m = m + 1$ and go back to Step \ref{step backward}.
	\State \label{step forward} Set $v^{\rm comp} = u_{m + 1} + g.$
 	\end{algorithmic}
\end{algorithm}

\begin{Remark}
The choices of $M$ and $\varepsilon$ in Step \ref{choice} are based on some
trial and error processes.
\end{Remark}

\begin{Remark}
Theorems \ref{thm 2.1} and \ref{thm 5.1} hold true even when the functions $%
v,$ $u$, $G$, $F$, $g_0$ and $g_1$ take complex values. By splitting %
\eqref{nonlinear eqn} and \eqref{nonlinear v} into real part and imaginary
part, we obtain a system of quasi-linear PDEs on the field of real numbers.
Then, we can apply the whole analysis for the case of a single equation to
the case of system of PDEs.
\end{Remark}

\section{A coefficient inverse problem in the frequency domain with back
scattering data in $\mathbb{R}^3$}

\label{sec CIP}

In this section, we introduce a method to solve the back scattering
inverse problem with multi-frequency data. This inverse problem that has
uncountable real-world applications. Solving a system like \eqref{nonlinear
eqn} plays an inportant role in our method.
Let $\Omega$ be the cube $(-R, R)^3 \subset \mathbb{R}^3$ where $R$ is a
positive number. Let $c \in C^1(\mathbb{R}^3)$ represent the dielectric
constant of $\mathbb{R}^3$. Assume that 
\begin{equation}
\left\{ 
\begin{array}{ll}
c(\mathbf{x}) = 1 & \mbox{if } \mathbf{x} \in \mathbb{R}^3 \setminus \Omega,
\\ 
c(\mathbf{x}) \geq 1 & \mbox{if } \mathbf{x} \in \Omega.%
\end{array}
\right.  \label{medium}
\end{equation}
Assumption \eqref{medium} can be understood as the dielectric constant of
the air (or vacuum) is scaled to be $1$. Let $[\underline k, \overline k]$
be an interval of wave number and let $\mathbf{x}_0 = (0, 0, -d)$, with $d >
R$, be a point located outside $\Omega$. Let $u = u(\mathbf{x}, k)$, $(%
\mathbf{x}, k) \in \mathbb{R}^3 \times [\underline k, \overline k]$,
represent the frequency-dependent wave. The function $u(\mathbf{x}, k)$ is
governed by the following problem 
\begin{equation}
\left\{ 
\begin{array}{rcll}
\Delta u(\mathbf{x}, k) + k^2 c(\mathbf{x}) u(\mathbf{x}, k) & = & -\delta(%
\mathbf{x} - \mathbf{x}_0) & \mathbf{x} \in \mathbb{R}^3, \\ 
\partial_{|\mathbf{x}|} u(\mathbf{x}, k) - \mathrm{i} k u(\mathbf{x}, k) & =
& o(|\mathbf{x}|^{-1}) & |\mathbf{x}| \to \infty%
\end{array}
\right.  \label{main eqn}
\end{equation}
where $\delta$ is the Dirac function. The partial differential differential
equation in \eqref{main eqn} is called the Helmholtz equation and the
asymptotic behavior of $u$ as $|\mathbf{x}| \to \infty$ is called the
Sommerfeld radiation condition. The Sommerfeld radiation condition
guarantees the existence and uniqueness of problem \eqref{main eqn}, see 
\cite[Chapter 8]{ColtonKress:2013}. We are interested in the following
problem.

\begin{problem}[Coefficient inverse problem from back scattering data]
\textit{\ Let  
\begin{equation}
\Gamma = \{\mathbf{x} = (x, y, -R): -R \leq x, y \leq R\} \subset \partial
\Omega  \label{gamma def}
\end{equation}
be the measurement site.  Given the measurements of  
\begin{equation}
f(\mathbf{x}, k) = u(\mathbf{x}, k) \quad \mbox{and} \quad g(\mathbf{x}, k)
= -\partial_{z} u(\mathbf{x}, k)  \label{5.3}
\end{equation}
for $(\mathbf{x}, k) \in \Gamma \times [\underline k, \overline k],$ 
determine the function $c(\mathbf{x}),$ $\mathbf{x} \in \Omega.$  \label{CIP}
}
\end{problem}

\begin{Remark}
    We refer the reader to our recent works \cite{Khoaelal:IPSE2021,KhoaKlibanovLoc:SIAMImaging2020,VoKlibanovNguyen:IP2020} in which we study a similar inverse problem in which the data are generated by a source moving along a straight line and only a single frequency was used. Unlike this, the data for the inverse problem under consideration, Problem \ref{CIP}, are generated by a single source and by multi-frequencies.
\end{Remark}

\begin{figure}[h!]
\begin{center}
\includegraphics[width = .5\textwidth]{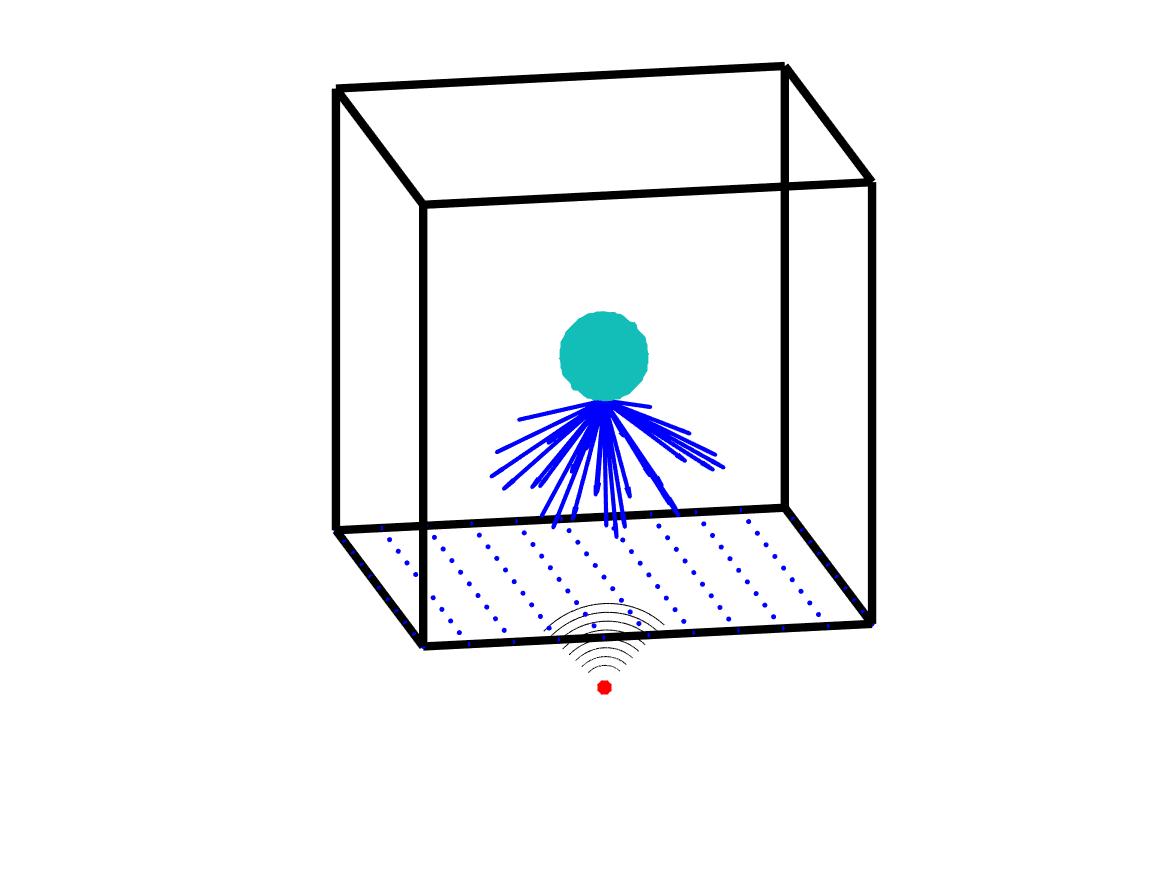}  {\tiny \ \put(-110, 95){$%
\Omega$}  \put(-56,42){$\Gamma$}  \put(-90,85){unknown object} 
\put(-66,62){back scattering wave}  \put(-80,24){source}  } 
\end{center}
\caption{A diagram for the experiment that leads to Problem \protect\ref{CIP}%
. The unknown object is located inside a box $\Omega$. An emitter (the red
dot), located outside $\Omega$, emits the incident wave. The incident wave
scatters when hitting the unknown object. The back scattering waves are
represented by blue arrows and collected on an array of detectors, located
on the part $\Gamma$ of $\partial \Omega$.}
\label{fig diagram}
\end{figure}
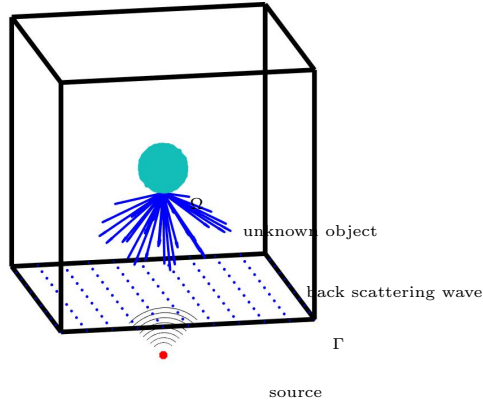

Problem \ref{CIP} arises from the following well-known experiment,
illustrated in Figure \ref{fig diagram}. Let an optical source illuminate
objects inside $\Omega$. The wave generated by the optical source is called
the incident wave. The incident wave hits the objects and scatters in all
directions. We collect the scattering wave on the part of $\partial \Omega$,
named as $\Gamma$, that receives the wave coming back from the objects.
Solving Problem \ref{CIP}, with these data, we obtain the dielectric
constant of the medium. This information is important in identifying the
objects. 
The forward problem corresponding to Problem \ref{CIP} is the problem of
computing the function $u(\mathbf{x}, k)$, $(\mathbf{x}, k) \in \Gamma
\times [\underline k, \overline k].$ To solve the forward problem, we first
model the incident wave by the point source 
\begin{equation}
u_0(\mathbf{x}, k) = \frac{\exp(\mathrm{i} k |\mathbf{x} - \mathbf{x}_0|)}{%
4\pi |\mathbf{x} - \mathbf{x}_0|} \quad (\mathbf{x}, k) \in \mathbb{R}^3
\times [\underline k, \overline k].
\end{equation}
It is well-known that 
\begin{equation}
\Delta u_0(\mathbf{x}, k) + k^2 u_0(\mathbf{x}, k) = -\delta(\mathbf{x} - 
\mathbf{x}_0) \quad \mathbf{x} \in \mathbb{R}^3, k \in [\underline k,
\overline k].  \label{2.4}
\end{equation}
Let 
\begin{equation}
u_{\mathrm{sc}}(\mathbf{x}, k) = u(\mathbf{x}, k) - u_0(\mathbf{x}, k) \quad 
\mathbf{x} \in \mathbb{R}^3, k \in [\underline k, \overline k]  \label{2.5}
\end{equation}
denote the scattering wave. It follows from \eqref{2.5} that $u(\mathbf{x},
k)$ is the sum of the scattering wave and the incident wave. We call the
function $u(\mathbf{x}, k)$ the total wave. Subtracting the differential
equation in \eqref{main eqn} from \eqref{2.4}, we obtain 
\begin{equation*}
\Delta u_{\mathrm{sc}}(\mathbf{x}, k) + k^2 u_{\mathrm{sc}}(\mathbf{x}, k) =
-k^2(c(\mathbf{x}) - 1) u(\mathbf{x}, k) \quad \mathbf{x} \in \mathbb{R}^3,
k \in [\underline k, \overline k].
\end{equation*}
Hence, see \cite[Chapter 8]{ColtonKress:2013}, 
\begin{equation}
u_{\mathrm{sc}}(\mathbf{x}, k) = k^2 \int_{\mathbb{R}^3} \frac{\exp(\mathrm{i%
} k |\mathbf{x} - \xi|)}{4\pi |\mathbf{x} - \xi|} (c(\xi) - 1) u(\xi, k)
d\xi \quad \mathbf{x} \in \mathbb{R}^3, k \in [\underline k, \overline k].
\label{2.6}
\end{equation}
Combining \eqref{2.5} and \eqref{2.6}, we arrive at the Lippmann-Schwinger
equation 
\begin{equation}
u(\mathbf{x}, k) = u_0(\mathbf{x}, k) + k^2 \int_{\mathbb{R}^3} \frac{\exp(%
\mathrm{i} k |\mathbf{x} - \xi|)}{4\pi |\mathbf{x} - \xi|} (c(\xi) - 1)
u(\xi, k) d\xi  \label{2.7}
\end{equation}
for all $\mathbf{x} \in \mathbb{R}^3, k \in [\underline k, \overline k].$ We
solve the integral equation \eqref{2.7} by the method in \cite%
{LechleiterNguyen:acm2014,Nguyen:anm2015}. In order to solve the inverse
problem, we need to impose the following condition.

\begin{Assumption}
\textit{\ The total wave $u(\mathbf{x}, k)$ is nonzero for all $\mathbf{x}
\in \Omega$ and $k \in [\underline k, \overline k].$  \label{assumption 2.1} 
}
\end{Assumption}

We provide here an example when Assumption \ref{assumption 2.1} holds true.
In this example, we assume that $c$ is in the class $C^{15}(\mathbb{R}^3)$.
Consider the Riemannian metric generated by $c$ 
\begin{equation}
d\tau = \sqrt{c(\mathbf{x})} d\mathbf{x}, \quad |d\mathbf{x}| = \sqrt{dx_1^2
+ dx_2^2 + dx_3^2}.  \label{Riemannian metric}
\end{equation}
Assume that for each point $\mathbf{x} \in \mathbb{R}^3$, the geodesic line
with respect to the Riemannian metric defined \eqref{Riemannian metric}
connecting $\mathbf{x}_0$ and $\mathbf{x}$ is unique, where $\mathbf{x}_0$
is the location of the emitter that emits the point source. Then, it was shown
in \cite{KlibanovRomanov:SIAMam2016} that 
\begin{equation}
u(\mathbf{x}, k) = A(\mathbf{x}) e^{\mathrm{i} k \tau(\mathbf{x})} + O\Big(%
\frac{1}{k}\Big) \quad \mbox{as } k \to \infty 
\label{atau}
\end{equation}
for $\mathbf{x} \in \Omega$ where $A$ is a function taking positive value
and $\tau$ is the travel time of the wave from $\mathbf{x}_0$ to $\mathbf{x}$%
. Hence, Assumption \ref{assumption 2.1} holds true when the wave number $k$
is sufficiently large. In the next section, we derive a system of nonlinear
partial differential equations. Solution of this system directly yields the
solution to Problem \ref{CIP}.

\section{A method to solve Problem \protect\ref{CIP}}

\label{sec 7}

Recall that $u = u(\mathbf{x}, k)$, $\mathbf{x} \in \Omega$, $k \in
[\underline k, \overline k]$ is the solution to \eqref{main eqn}. Assume
that Assumption \ref{assumption 2.1} holds true. Define 
\begin{equation}
v(\mathbf{x}, k) = \frac{1}{k^2} \log \frac{u(\mathbf{x}, k)}{u_0(\mathbf{x}%
, k)} \quad \mbox{for all } \mathbf{x} \in \Omega, k \in [\underline k,
\overline k].  \label{fn v}
\end{equation}
%\begin{Remark}
Although $u/u_0$ takes complex values, the function $v$ can be defined.
Employing \eqref{atau} and  assuming that both $\underline k$ and $\overline k$ are large, we define the function $v$ as
\begin{align*}
    v(\x, k) 
    &= \frac{1}{k^2} \big[\log u(\x, k) - \log u_0(\x, k)\big]
    \\
    &= \frac{1}{k^2}
    \Big[
        \ln A(\x) - \ln \frac{1}{4\pi|\x - \x_0|} + \ik (\tau(\x) - |\x - \x_0|)
        + O(1/k)
    \Big]
\end{align*}
for all $\x \in \Omega, k \in [\underline k, \overline k].$ 
%\end{Remark}
We now derive a differential equation for the function $v$. It follows from %
\eqref{fn v} that 
\begin{equation}
\nabla v(\mathbf{x}, k) = \frac{1}{k^2} \Big(\frac{\nabla u(\mathbf{x}, k)}{%
u(\mathbf{x}, k)} - \frac{\nabla u_0(\mathbf{x}, k)}{u_0(\mathbf{x}, k)}%
\Big) \quad \mbox{for all } \mathbf{x} \in \Omega, k \in [\underline k,
\overline k].  \label{3.5}
\end{equation}
Taking the divergence of \eqref{3.5} gives 
\begin{equation}
\Delta v(\mathbf{x}, k) = \frac{1}{k^2} \Big( \frac{\Delta u(\mathbf{x}, k)}{%
u(\mathbf{x}, k)} - \Big(\frac{\nabla u(\mathbf{x}, k)}{u(\mathbf{x}, k)}%
\Big)^2 -\frac{\Delta u_0(\mathbf{x}, k)}{u_0(\mathbf{x}, k)} + \Big(\frac{%
\nabla u_0(\mathbf{x}, k)}{u_0(\mathbf{x}, k)}\Big)^2 \Big)  \label{3.6}
\end{equation}
for all $\mathbf{x} \in \Omega$, $k \in [\underline k, \overline k]$. Since $%
u(\mathbf{x}, k)$ satisfies the differential equation in \eqref{main eqn}
and $u_0(\mathbf{x}, k)$ satisfies the differential equation in \eqref{main
eqn} with $c(\mathbf{x})$ replaced by $1$, we have 
\begin{equation}
\frac{\Delta u(\mathbf{x}, k)}{u(\mathbf{x}, k)} = -k^2 c(\mathbf{x}) \quad %
\mbox{and} \quad \frac{\Delta u_0(\mathbf{x}, k)}{u_0(\mathbf{x}, k)} = -k^2
\label{3.7}
\end{equation}
for all $\mathbf{x} \in \Omega$, $k \in [\underline k, \overline k]$. It
follows from \eqref{3.5}, \eqref{3.6} and \eqref{3.7} that 
\begin{align*}
\Delta v(\mathbf{x}, k) &= \frac{1}{k^2} \Big[ -k^2(c(\mathbf{x}) - 1) - %
\Big( k^2 \nabla v(\mathbf{x}, k) + \frac{\nabla u_0(\mathbf{x}, k)}{u_0(%
\mathbf{x}, k)} \Big)^2 + \Big(\frac{\nabla u_0(\mathbf{x}, k)}{u_0(\mathbf{x%
}, k)}\Big)^2\Big] \\
&= -(c(\mathbf{x}) - 1) - k^2 (\nabla v(\mathbf{x}, k))^2 - \frac{2\nabla v(%
\mathbf{x}, k) \cdot \nabla u_0(\mathbf{x}, k)}{u_0(\mathbf{x}, k)}
\end{align*}
for all $\mathbf{x} \in \Omega$, $k \in [\underline k, \overline k]$. We
obtain 
\begin{equation}
\Delta v(\mathbf{x}, k) + k^2(\nabla v(\mathbf{x}, k))^2 + \frac{2\nabla v(%
\mathbf{x}, k) \cdot \nabla u_0(\mathbf{x}, k)}{u_0(\mathbf{x}, k)} = -c(%
\mathbf{x}) + 1  \label{6.5}
\end{equation}
for all $\mathbf{x} \in \Omega$, $k \in [\underline k, \overline k]$.
Differentiate \eqref{6.5} with respect to $k$. We have 
\begin{multline}
\Delta \partial_k v(\mathbf{x}, k) + 2 k^2 \nabla v(\mathbf{x}, k) \cdot
\nabla \partial_k v(\mathbf{x}, k) + 2k(\nabla v(\mathbf{x}, k))^2 \\
+ 2\nabla \partial_k v(\mathbf{x}, k) \cdot \frac{\nabla u_0(\mathbf{x}, k)}{%
u_0(\mathbf{x}, k)} + 2\nabla v(\mathbf{x}, k) \partial_k \frac{\nabla u_0(%
\mathbf{x}, k)}{u_0(\mathbf{x}, k)} = 0  \label{6.6}
\end{multline}
for all $\mathbf{x} \in \Omega$, $k \in [\underline k, \overline k]$.

Let $\{\Psi_m\}_{m \geq 1}$ be the orthonormal basis of $L^2(\Omega)$
defined in \cite{Klibanov:jiip2017}. This basis is constructed as follows.
For each $m \geq 1,$ define the function $\phi_m(k) = k^{m - 1}e^{k -
(\overline k + \underline k)/2}$ for all $k \in [\underline k, \overline k].$
It is clear that the set $\{\phi_m\}_{m \geq 1}$ is complete in $[\underline
k, \overline k].$ We then apply the Gram-Schmidt orthonormalization process
to this set to obtain the basis $\{\Psi_m\}_{m \geq 1}$.

We derive an approximation model for the solution $v$ to \eqref{6.6} as
follows. For each $\mathbf{x} \in \Omega$ and $k \in [\underline k,
\overline k],$ we write 
\begin{equation}
v(\mathbf{x}, k) = \sum_{i = 1}^\infty v_i(\mathbf{x}) \Psi_i(k) \simeq
\sum_{i = 1}^N v_i(\mathbf{x}) \Psi_i(k)  \label{6.7}
\end{equation}
for some cut-off number $N$, determined later in section \ref{sec num},
where 
\begin{equation}
v_i(\mathbf{x}) = \int_{\underline k}^{\overline k} v(\mathbf{x}, \kappa)
\Psi_i(\kappa) d\kappa \quad i \in \{1, \dots, N\}.  \label{6.8}
\end{equation}
In this approximation context, 
\begin{equation}
v_k(\mathbf{x}, k) = \sum_{i = 1}^N v_i(\mathbf{x}) \Psi_i^{\prime }(k)
\quad \mbox{for all } \mathbf{x} \in \Omega, k \in [\underline k, \overline
k].  \label{6.9}
\end{equation}
Plugging \eqref{6.7} and \eqref{6.9} into \eqref{6.6} gives 
\begin{multline}
\sum_{i = 1}^N \Delta v_{i}(\mathbf{x}) \Psi_i^{\prime }(k) + 2\sum_{i, j =
1}^N \nabla v_i(\mathbf{x}) \cdot \nabla v_j(\mathbf{x}) \Big(k^2
\Psi_i(k)\Psi_j^{\prime }(k) + k \Psi_i(k)\Psi_j(k)\Big) \\
+ 2\sum_{i = 1}^N \nabla v_i(\mathbf{x}) \cdot \Big( \Psi_i^{\prime }(k)%
\frac{\nabla u_0(\mathbf{x}, k)}{u_0(\mathbf{x}, k)} + \Psi_i(k) \partial_k 
\frac{\nabla u_0(\mathbf{x}, k)}{u_0(\mathbf{x}, k)} \Big) = 0  \label{6.10}
\end{multline}
for all $\mathbf{x} \in \Omega.$ For each $l \in \{1, \dots, N\},$
multiplying $\Psi_l(k)$ to both sides of \eqref{6.10}, we have 
\begin{equation}
\sum_{i = 1}^N s_{li} \Delta v_i(\mathbf{x}) + \sum_{i, j = 1}^N
a_{lij}\nabla v_i(\mathbf{x}) \cdot \nabla v_j(\mathbf{x}) + \sum_{i = 1}^N
B_{li}(\mathbf{x}) \cdot \nabla v_i(\mathbf{x}) = 0  \label{6.11}
\end{equation}
where 
\begin{equation}
\left\{  
\begin{array}{l}
s_{li} = \displaystyle\int_{\underline k}^{\overline k} \Psi_i^{\prime }(k)
\Psi_l(k) dk, \\ 
a_{lij} = \displaystyle 2 \int_{\underline k}^{\overline k} \Big(k^2
\Psi_i(k)\Psi_j^{\prime }(k) + k \Psi_i(k)\Psi_j(k)\Big)\Psi_l(k)dk, \\ 
B_{li}(\mathbf{x}) = \displaystyle 2 \int_{\underline k}^{\overline k}\Big( %
\Psi_i^{\prime }(k)\frac{\nabla u_0(\mathbf{x}, k)}{u_0(\mathbf{x}, k)} +
\Psi_i(k) \partial_k \frac{\nabla u_0(\mathbf{x}, k)}{u_0(\mathbf{x}, k)} %
\Big)\Psi_l(k)dk%
\end{array}
\right. 
\label{sab}
\end{equation}
for all $i, j, l \in \{1, \dots, N\}$ and $\mathbf{x} \in \Omega$.

We next compute the boundary information for $V$. For all $\mathbf{x} \in
\Gamma$ and $k \in [\underline k, \overline k]$, since $u(\mathbf{x}, k) = f(%
\mathbf{x}, k)$ where $f$ is the data for the inverse problem under
consideration (see \eqref{5.3}), it follows from \eqref{fn v} that $v(%
\mathbf{x}, k) = \frac{1}{k^2} \log \frac{f(\mathbf{x}, k)}{u_0(\mathbf{x},
k)}$. Hence, by \eqref{6.8}, we have 
\begin{equation*}
v_i(\mathbf{x}) = \int_{\underline k}^{\overline k} v(\mathbf{x}, \kappa)
\Psi_i(\kappa) d\kappa = \int_{\underline k}^{\overline k} \frac{1}{\kappa^2}
\log \frac{f(\mathbf{x}, \kappa)}{u_0(\mathbf{x}, \kappa)} \Psi_i(\kappa)
d\kappa \quad i \in \{1, \dots, N\}
\end{equation*}
for all $\mathbf{x} \in \Gamma.$ Since we only measure the wave $u$ on $%
\Gamma,$ we complement $v_i(\mathbf{x}) = 0$ for all $%
i \in \{1, \cdots, N\}$ and $\mathbf{x} \in \partial \Omega \setminus \Gamma.
$ The boundary value of $V_N = (v_1, v_2, \dots, v_N)^{\rm T}$ is given by 
\begin{equation}
V_N(\x) = g_1(\mathbf{x}) = \left\{ 
\begin{array}{ll}
\begin{array}{c}
\Big[
\displaystyle \int_{\underline k}^{\overline k} \frac{1}{\kappa^2} \log 
\frac{f(\mathbf{x}, \kappa)}{u_0(\mathbf{x}, \kappa)} \Psi_i(\kappa) d\kappa%
\Big]_{i = 1}^N
\end{array}
 & \mathbf{x} \in \Gamma, \\ 
0 & \mathbf{x} \in \partial{\Omega} \setminus \Gamma%
\end{array}
\right.  \label{6.12}
\end{equation}
On the other hand, by \eqref{6.8}, for all $\mathbf{x} \in \Gamma$ and $k
\in [\underline k, \overline k]$, 
\begin{equation*}
\partial_{\nu} v(\mathbf{x}, k) = \frac{1}{k^2} \Big(\frac{\partial_{\nu} u(%
\mathbf{x}, k)}{u(\mathbf{x}, k)} - \frac{\partial_{\nu} u_0(\mathbf{x}, k)}{%
u_0(\mathbf{x}, k)}\Big) = \frac{1}{k^2} \Big(\frac{g(\mathbf{x}, k)}{f(%
\mathbf{x}, k)} - \frac{\partial_{\nu} u_0(\mathbf{x}, k)}{u_0(\mathbf{x}, k)%
}\Big). 
\end{equation*}
Therefore, by \eqref{6.8}, 
\begin{align}
\partial_{\nu} v_i(\mathbf{x}) 
&=\int_{\underline k}^{\overline k}
\partial_{\nu}v(\mathbf{x}, \kappa) \Psi_i(\kappa) d\kappa  \nonumber
\\
&=
\int_{\underline k}^{\overline k} \frac{1}{\kappa^2} \Big(\frac{g(\mathbf{x}%
, \kappa)}{f(\mathbf{x}, \kappa)} - \frac{\partial_{\nu} u_0(\mathbf{x},
\kappa)}{u_0(\mathbf{x}, \kappa)}\Big) \Psi_i(\kappa)d\kappa
\end{align}
for all $i \in \{1, \dots, N\}$ and $\mathbf{x} \in \Gamma.$ 
By setting 
\begin{equation}
g_0(\mathbf{x}) = \left[  \int_{\underline k}^{\overline k} \frac{1}{\kappa^2%
} \Big(\frac{g(\mathbf{x}, \kappa)}{f(\mathbf{x}, \kappa)} - \frac{%
\partial_{\nu} u_0(\mathbf{x}, \kappa)}{u_0(\mathbf{x}, \kappa)}\Big) %
\Psi_i(\kappa) d\kappa  \right]_{i = 1}^N, 
\label{g1i}
\end{equation}
we have 
\begin{equation}
\partial_{\nu} V(\mathbf{x}) = g_0(\mathbf{x}) \quad \mbox{for all } \mathbf{%
x} \in \Gamma.  \label{6.14}
\end{equation}
In summary, the vector $V_N = (v_1, v_2, \dots, v_N)^{\rm T}$ satisfies a Cauchy like boundary value problem
\begin{equation}
    \left\{
        \begin{array}{rcll}
             \displaystyle\sum_{i = 1}^N s_{li} \Delta v_i(\mathbf{x}) + \sum_{i, j = 1}^N
a_{lij}\nabla v_i(\mathbf{x}) \cdot \nabla v_j(\mathbf{x})\\
\hspace{3cm}+
\displaystyle\sum_{i = 1}^N  
B_{li}(\mathbf{x}) \cdot \nabla v_i(\mathbf{x}) &=& 0,
             &\quad \x \in \Omega,\\
             V_N(\x)  &=& g_1(\x)  &\quad \x \in \partial \Omega,
             \\
             \partial_{\nu} V_N(\x)   &=& g_0(\x) &\quad \x \in \Gamma.
        \end{array}
    \right.
    \label{quasi-system}
\end{equation}
% \begin{multline}
%     \sum_{i = 1}^N s_{li} \Delta v_i(\mathbf{x}) + \sum_{i, j = 1}^N
% a_{lij}\nabla v_i(\mathbf{x}) \cdot \nabla v_j(\mathbf{x})
% \sum_{i = 1}^N
% \\
% +
% B_{li}(\mathbf{x}) \cdot \nabla v_i(\mathbf{x}) = 0
% \quad \mbox{for all } \x \in \Omega,
% \end{multline}
% and
% \begin{equation}
% \left\{
%     \begin{array}{rcll}
%       V_N(\x)  &=& g_1(\x)  & \x \in \partial \Omega, \\
%       \partial_{\nu} V_N(\x)   &=& g_0(\x) &\x \in \Gamma.
%     \end{array}
% \right.
% \label{boundarysystem}
% \end{equation}
where $s_{ij}$, $a_{lij}$, $B_{li}$, $1 \leq i, j, l \leq N$ are given in \eqref{sab},  $g_1$ and $g_0$, are respectively defined in \eqref{6.12}  and \eqref{g1i}.

\begin{Remark}
Let $(\tilde s_{li})_{l, i = 1}^N$ denote $S^{-1}$. Problem \eqref{quasi-system} can be rewritten as a particular form of \eqref{nonlinear eqn} as
\begin{equation*}
    \left\{
        \begin{array}{rcll}
             \displaystyle\sum_{i = 1}^N  \Delta v_i(\mathbf{x}) + \sum_{i, j = 1}^N
\tilde
s_{li} a_{lij}\nabla v_i(\mathbf{x}) \cdot \nabla v_j(\mathbf{x})\\
\hspace{3cm}+
\displaystyle\sum_{i = 1}^N  
\tilde
s_{li} B_{li}(\mathbf{x}) \cdot \nabla v_i(\mathbf{x}) &=& 0,
             &\quad \x \in \Omega,\\
             V_N(\x)  &=& g_1(\x)  &\quad \x \in \partial \Omega,
             \\
             \partial_{\nu} V_N(\x)   &=& g_0(\x) &\quad \x \in \Gamma.
        \end{array}
    \right.
\end{equation*}
% \begin{equation*}
% \sum_{i = 1}^N \Delta v_i(\mathbf{x}) + \sum_{i, j = 1}^N \tilde
% s_{li}a_{lij}\nabla v_i(\mathbf{x}) \cdot \nabla v_j(\mathbf{x}) + \sum_{i =
% 1}^N \tilde s_{li} B_{li}(\mathbf{x}) \cdot \nabla v_i(\mathbf{x}) = 0
% \end{equation*}
%for $\mathbf{x} \in \Omega$ and the boundary conditions \eqref{boundarysystem}. 
However, we numerically observe that solving \eqref{quasi-system} provides better numerical
solutions.
\end{Remark}

%
%We note that quasi-linear equations \eqref{6.11} with the mixed boundary data \eqref{6.12} and \eqref{6.14} is a particular case of \eqref{nonlinear eqn}. 
%To solve it, we employ Algorithm \ref{alg 1}. 
%Let $g$ be a vector that satisfies \eqref{6.12} and \eqref{6.14} (we will present a method to compute $g$ in section \ref{subsec initial}). 
%Define $U = (u_1, \dots, u_N) = V - g$ as in \eqref{change}. The system \eqref{6.11} becomes
%\begin{equation}
%	S\Delta U =  G(\x, U(\x), \nabla U(\x)) \quad \x \in \Omega
%\end{equation}
%where
%\begin{equation}
%	G(\x, U(\x), \nabla U(\x)) = S\Delta g + F(\x, U(\x) + g(\x), \nabla (U(\x) + g(\x))).
%\end{equation}

Let $\mu_\lambda(x, y, z) = e^{\lambda (z + r)^2}$ for a number $r > 1$ be a
Carleman weight function. We refer the reader to \cite[Theorem 3.1]%
{KhoaKlibanovLoc:SIAMImaging2020} for the proof of Carleman estimate %
\eqref{Car est} with this Carleman weight function in 3D. Thus, we can find
the solution to the system of quasi-linear equations \eqref{quasi-system} by the convexification method,
see Algorithm \ref{alg 1}.

\section{Numerical tests}\label{sec num}

In this section, we numerically study Problem \ref{p1} and Problem \ref{CIP}.

\subsection{Numerical study for Problem \ref{p1}} 

We present an example in which we apply  convexification method to compute the solution to problem of the form \eqref{nonlinear eqn}.
For simplicity, we consider the case $n = 2$ and
the computational domain $\Omega$ is chosen to be $(-1, 1)^2$. 
We choose the set $\Gamma$ to be $\{(x, y = -1): |x| \leq 1\} \subset \partial \Omega$, on which we impose the Neumann boundary condition for the solution $v$.
We numerically test the convexification method in the finite difference scheme.
That means we compute the values of the solution  on the grid
\[
    \mathcal{G} = \big\{
        (x_i = -1 + (i-1)d_{\x}, y_j = -1 + (j - 1)d_{\x}):
        1 \leq i, j \leq N_{\x}
    \big\}
\] where $d_\x = \frac{2}{N_{\x}-1}$ and $N_{\x}$ is an integer. In our numerical study, $N_{\x} = 41.$

\begin{figure}[!ht]
    \subfloat[The function $v^*$]{\includegraphics[width=.3\textwidth]{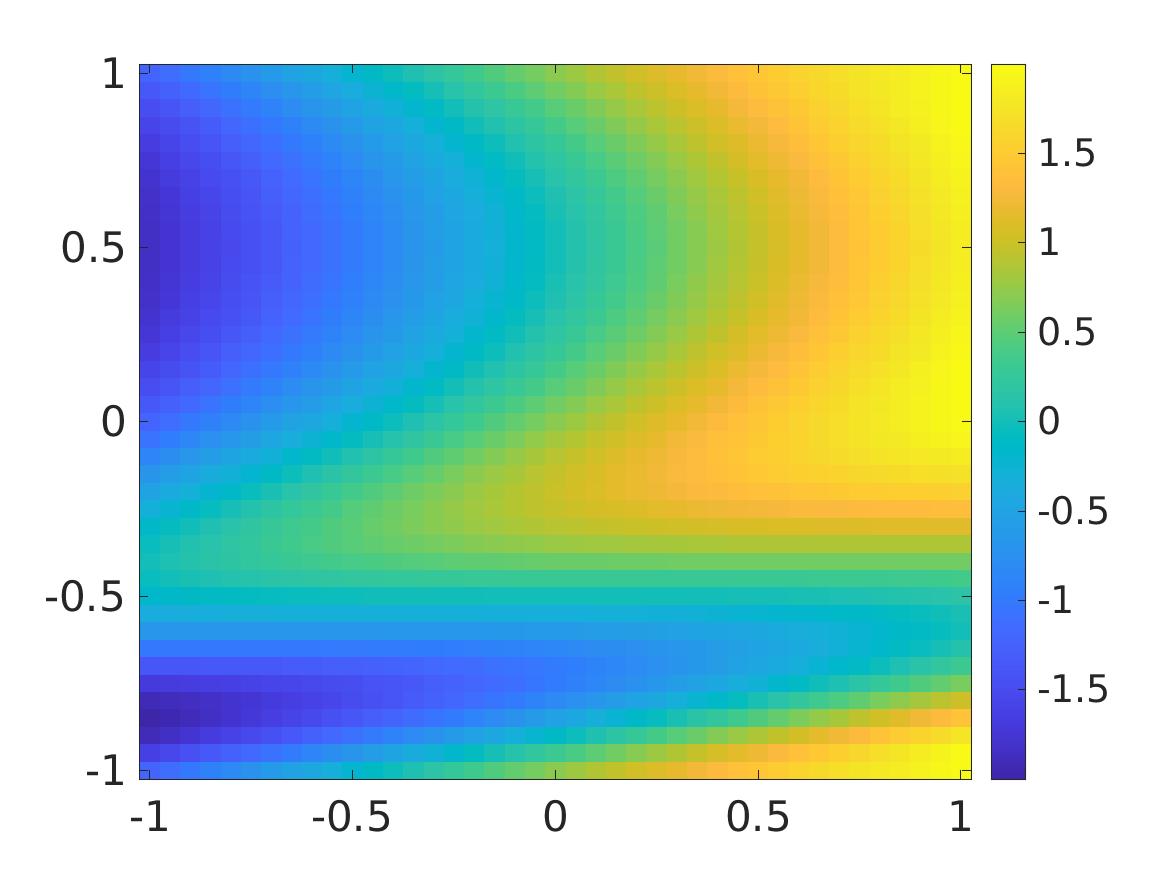}}
   
    \subfloat[The function $v^{\rm comp}$. $\delta = 5\%$.]{\includegraphics[width=.3\textwidth]{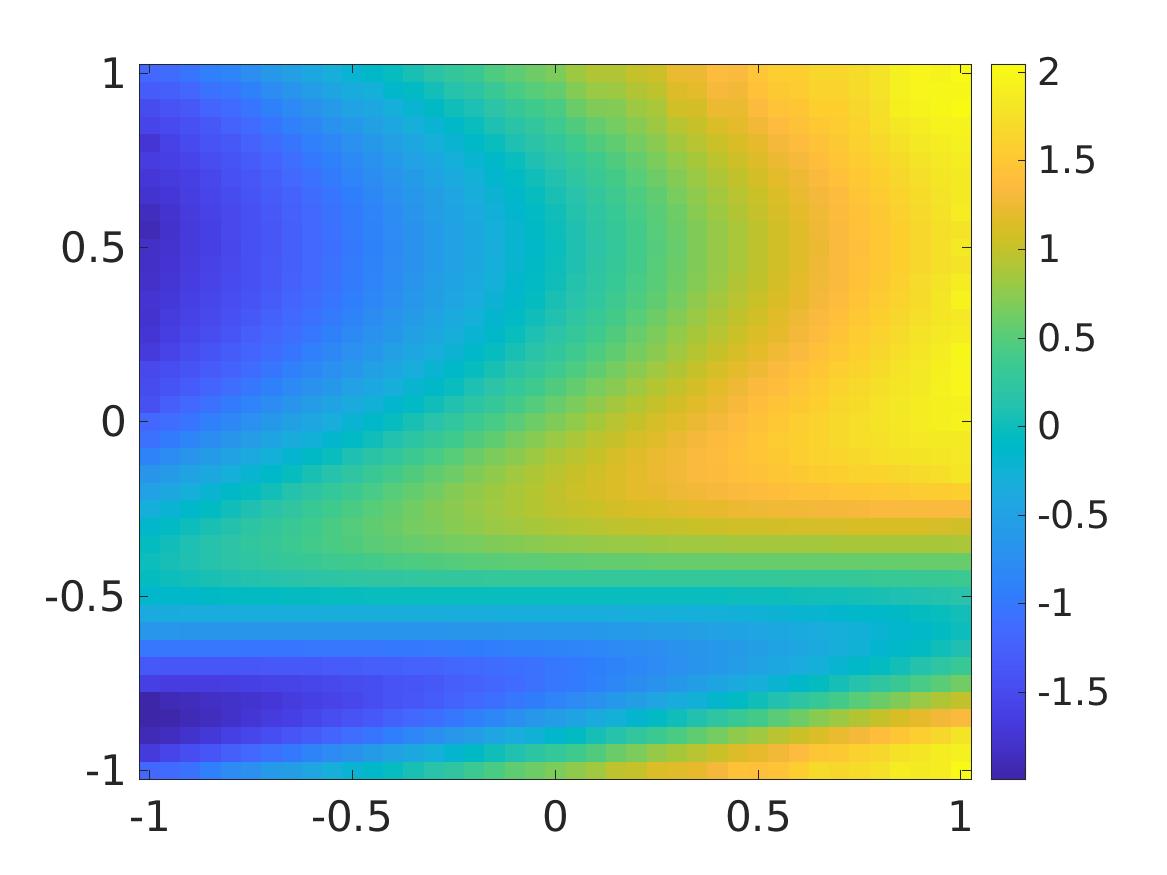}}
    \quad
    \subfloat[The function $v^{\rm comp}$. $\delta = 10\%$.]{\includegraphics[width=.3\textwidth]{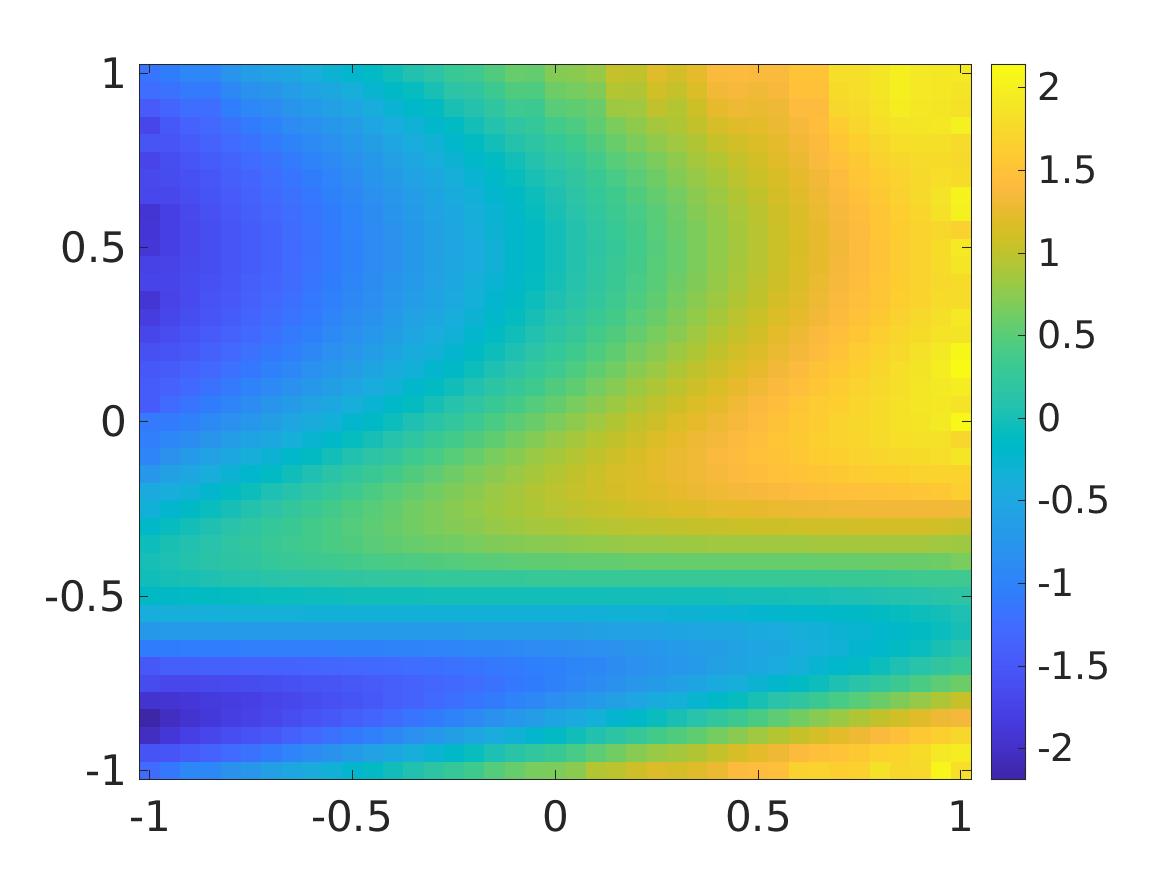}}
    \quad
    \subfloat[The function $v^{\rm comp}$. $\delta = 20\%$.]{\includegraphics[width=.3\textwidth]{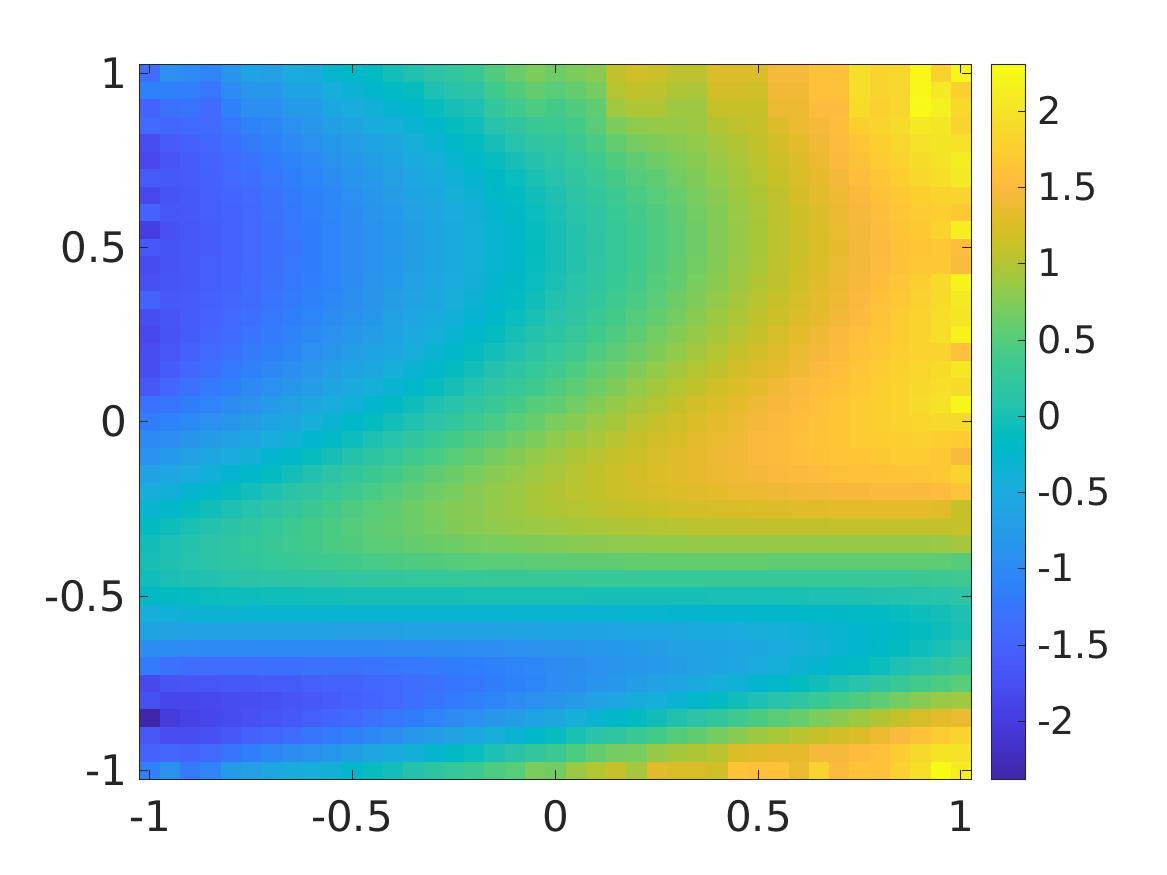}}
   
    \subfloat[The relative difference $\frac{|v^* - v^{\rm comp}|}{\|v^*\|_{L^{\infty}(\Omega)}}$. $\delta = 5\%$]{\includegraphics[width=.3\textwidth]{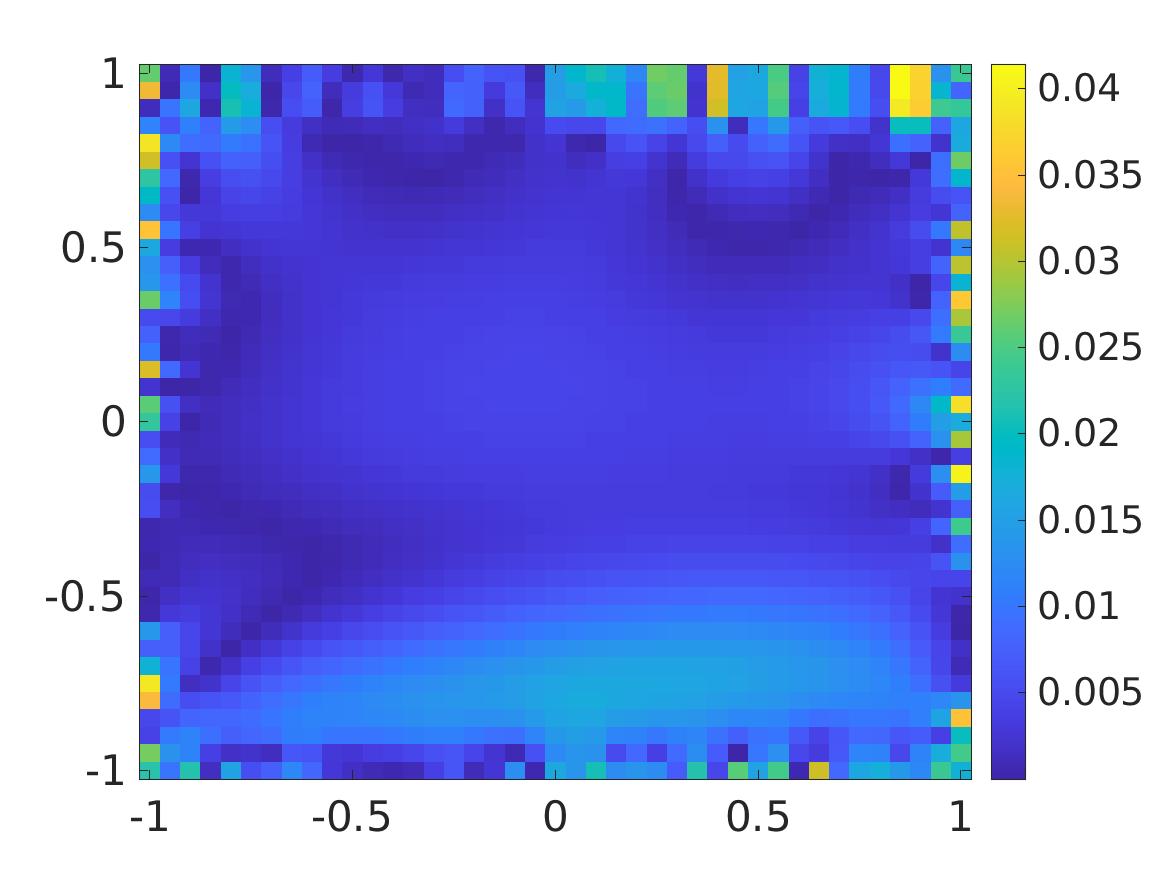}}
    \quad
    \subfloat[The relative difference $\frac{|v^* - v^{\rm comp}|}{\|v^*\|_{L^{\infty}(\Omega)}}$. $\delta = 10\%$]{\includegraphics[width=.3\textwidth]{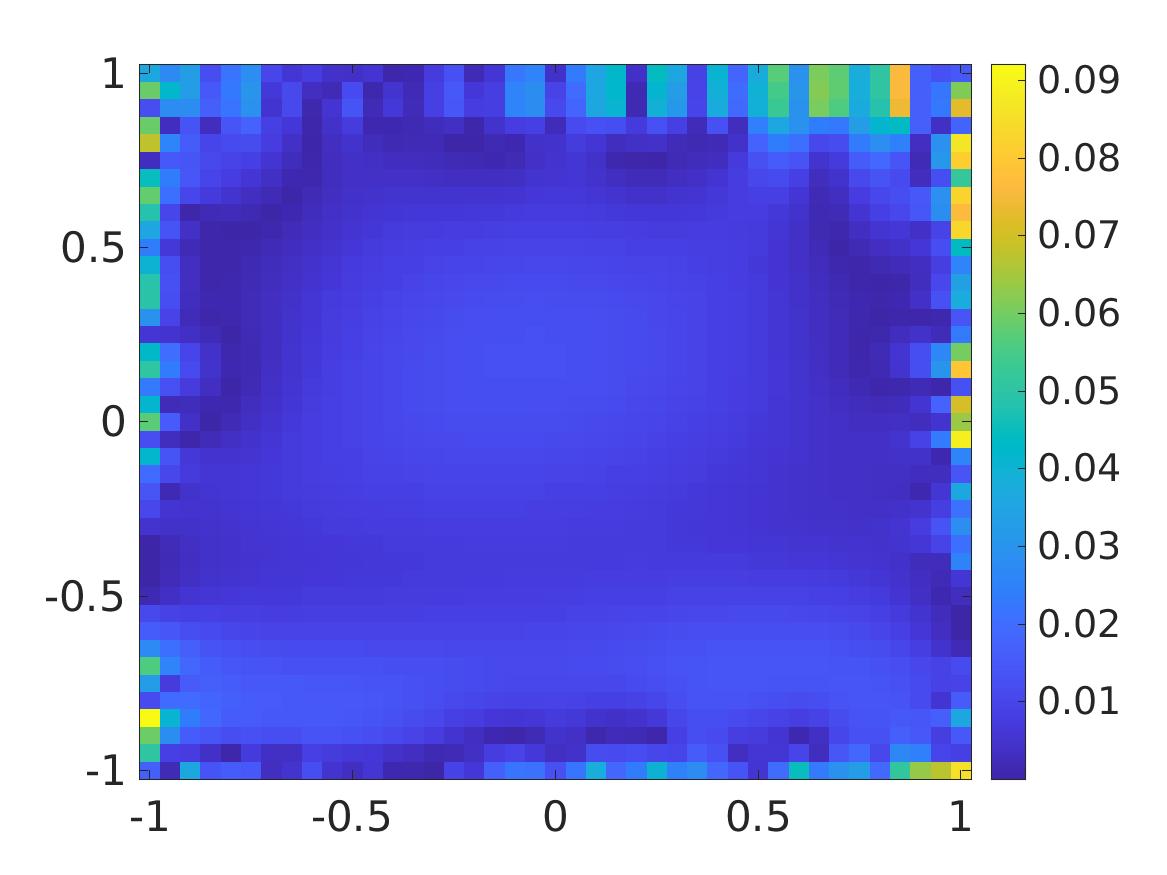}}
    \quad
    \subfloat[The relative difference $\frac{|v^* - v^{\rm comp}|}{\|v^*\|_{L^{\infty}(\Omega)}}$. $\delta = 20\%$]{\includegraphics[width=.3\textwidth]{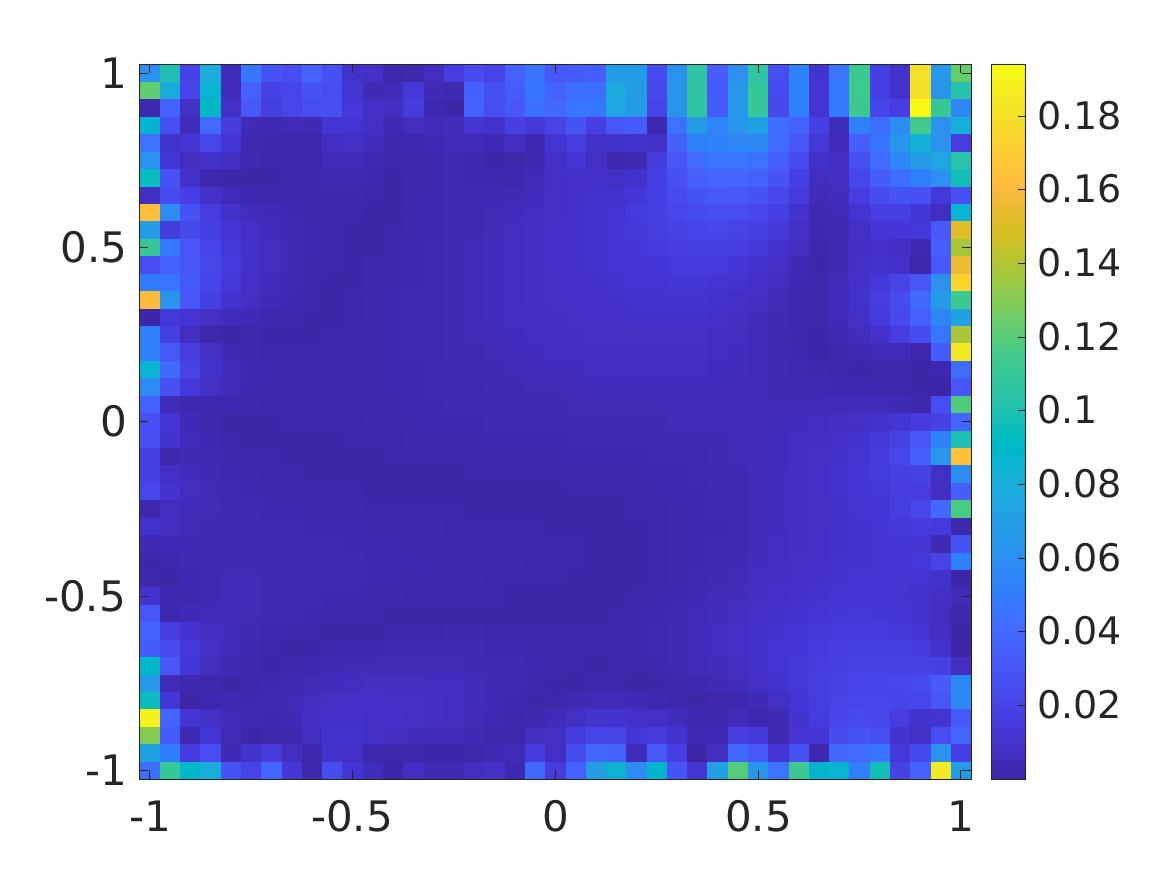}}
    \caption{Solutions to \eqref{nonlinear eqn} when $G$ and $f^*$ and $g^*$ are given in \eqref{G3} and \eqref{fg3} respectively.}
    \label{fig7_1}
\end{figure}

In computation, we use the Carleman weight function $ e^{-\lambda(R + 1.5)^2}e^{\lambda (y - 1.5)^2}$ where $\lambda = 1.1$.
This Carleman weight function is the 2D version of the one used in Section \ref{sec num2}. The regularization term is chosen to be $\epsilon = 10^{-6}.$
The details in implementation including the discretizing  the cost functional $J_{\epsilon, \lambda}$ and the choice of the initial guess are similar to the ones in Section \ref{sec num2}. We do not present in details here.  
To illustrate the efficiency of Algorithm \ref{alg 1}, we compute solution to \eqref{nonlinear eqn} when the nonlineariry $G$ is given by
\begin{multline}
   G(\x, s, p) =  - \sqrt{|p|^2 + 1} -\sin \left( x+\pi\, \left( y- 0.5 \right) ^{2} \right) +2\,\pi\,\cos
 \left( x+\pi\, \left( y- 0.5 \right) ^{2} \right)
 \\
 -4\,{\pi}^{2}
 \left( y- 0.5 \right) ^{2}\sin \left( x+\pi\, \left( y- 0.5 \right) ^
{2} \right) 
+\Big[  \left( 1+\cos \left( x+\pi\, \left( y- 0.5
 \right) ^{2} \right)  \right) ^{2}
 \\
 +4\,{\pi}^{2} \left( y- 0.5
 \right) ^{2} \left( \cos \left( x+\pi\, \left( y- 0.5 \right) ^{2}
 \right)  \right) ^{2}+1\Big]^{1/2}
 \label{G3}
\end{multline}
for all $\x = (x, y) \in \Omega$, $s \in \R$ and $p \in \R^d$. 
The exact boundary data are given by
\begin{equation}
    \left\{
        \begin{array}{ll}
             g_1^*(x, y) = x + \sin\big(x + \pi(y - 0.5)^2\big) &(x, y) \in \partial \Omega,\\
             g_0^*(x, y) = -2\pi (y - 0.5) \cos\big(x + \pi(y - 0.5)^2\big)& (x, y) \in \Gamma.
        \end{array}
    \right.
   \label{fg3}
\end{equation}
We add noise into the boundary data by the following formulas
\begin{equation*}
    g_1 = g_1^*(1 + \delta {\rm rand}) \quad \mbox{and }
    \quad g_0 = g_0^*(1 + \delta {\rm rand})
  %  \label{54}
\end{equation*}
where rand is a function taking uniformly distributed random numbers in the range $[-1, 1].$ 
The noise level $\delta$ is set to be $5\%$, $10\%$ and $20\%$.
The exact solution to \eqref{nonlinear eqn} in this test is $v^*(x, y) = x + \sin\big(x + \pi(y - 0.5)^2\big)$ for all $\x = (x, y) \in \Omega.$
The error in computation is given in  Table \ref{tab:my_label}.
\begin{table}[!ht]
    \centering
    \begin{tabular}{|c|c|}
    \hline
         $\delta$& relative error $\frac{\|v^* - v^{\rm comp}\|_{L^{\infty}(\Omega)}}{\|v^*\|_{L^{\infty}(\Omega)}}$ \\
         \hline
         5\% & 4.14\%\\
         \hline
         10\% & 9.21\%\\
         \hline
         20\% &19.40\%\\
         \hline
    \end{tabular}
    \caption{The relative error in our computation agaist the noise level $\delta$ contained in the boundary data.}
    \label{tab:my_label}
\end{table}
The graphs of the exact solution $v^*$, computed solution $v^{\rm comp}$ and their relative differences $\frac{|v^* - v^{\rm comp}|}{\|v^*\|_{L^{\infty}(\Omega)}},$ when $\delta = 5\%,$ $10\%$ and $20\%$ are displayed in Figure \ref{fig7_1}.

It is evident from Table \ref{tab:my_label} and Figure \ref{fig7_1} that the convexification method  delivers reliable solutions to quasi-linear elliptic equations. The errors in computation are compatible with the noise level and they occur on $\partial \Omega$ where the noise takes place.

%To minimize $J_{\epsilon, \lambda}$, we use the optimization package of Matlab to minimize the cost functional $J_{\epsilon, \lambda}$ by the gradient descent method. We perform two tests.

% \noindent{\bf Test 1.} We solve the problem
% \begin{multline*}
%     \Delta v + 10 v + v_x^2 - v_y^2 
%     =-{\pi}^{2}\sin \left( \pi\,x+y \right) +9\,\cos \left( \pi\,y-x
%  \right)
%  \\ +9\,\sin \left( \pi\,x+y \right) 
%  -{\pi}^{2}\cos \left( \pi\,y
% -x \right) + \left( \pi\,\cos \left( \pi\,x+y \right) +\sin \left( \pi
% \,y-x \right)  \right) ^{2}
% \\- \left( \cos \left( \pi\,x+y \right) -\pi
% \,\sin \left( \pi\,y-x \right)  \right) ^{2}
% \end{multline*}
% for all $x, y \in \Omega.$
% The boundary conditions are given by
% \[
%     v(x, y) = \sin(\pi x + y) + \cos(x - \pi y)
% \]
% for all $(x, y) \in \partial \Omega$ and 
% \[
%     \partial_{\nu} v(x, y) = \cos \left( \pi\,x+y \right) -\pi\,\sin \left( \pi\,y-x \right) 
% \]
% for all $(x, y) \in \Gamma$.

\subsection{Numerical study for Problem \ref{CIP}} \label{sec num2}
In this section, we present some numerical solutions to Problem \ref{CIP}.
The numerical examples we present in this section illustrate the efficiency
of the gradient descent method for the convexification described in section %
\ref{sec convex}. 
%Moreover, these examples serve as the numerical proof of the convergence of the gradient descent method in section \ref{sec gradient des}.
Especially, we will show that the presence of the Carleman weight function
in the objective function is crucial. That means, without involving the
Carleman estimate, the descent gradient method does not deliver good
numerical solutions to the problem of minimizing our nonconvex objective
functional.

\subsubsection{The forward problem}

The experimental setting for Problem \ref{CIP} is as follows. Let $R = 1$
and $\Omega = (-R, R)^3$. The source location is placed at $(0, 0, -4)$. The
interval of wavenumbers is $[\pi, 2\pi]$, which correspond to the interval
of wavelengths $[0.5, 1].$ In order to generate the simulated data, we use
the finite difference method in which we decompose $\Omega$ as the uniform
partition with grid points 
\begin{multline*}
\mathcal{G} = \{(x_i = -R + (i - 1)h, y_j = -R + (j - 1)h, z_l = -R + (l -
1)h):
\\1 \leq i, j, k \leq N_\mathbf{x}\} 
\end{multline*}
where $N_\mathbf{x} = 41$ and $h = 2R/(N_\mathbf{x} - 1)$. We also split the
interval of wavenumbers to the uniform partition 
\begin{equation*}
\mathcal{K} = \{k_1 = \underline k, k_2, \dots, k_{N_k} = \overline k\} 
\end{equation*}
where $k_i = \underline k + (i - 1)(\overline k - \underline k)/(N_k - 1)$
and $N_k = 121.$ The forward problem is solved via solving the
Lippmann-Schwinger equation \eqref{2.7} by the method in \cite%
{LechleiterNguyen:acm2014,Nguyen:anm2015}. Denote by $u^*(\mathbf{x}, k)$, $%
\mathbf{x} \in \mathcal{G}$, $k \in \mathcal{K}$ the obtained solution. The
noisy data for Problem \ref{CIP} is given by 
\begin{equation*}
f(\mathbf{x}, k) = u(\mathbf{x}, k)(1 + \delta \mathrm{rand}) \quad g(%
\mathbf{x}, k) = -\partial_z u^*(\mathbf{x}, k)(1 + \delta \mathrm{rand}) 
\end{equation*}
for all $\mathbf{x} \in \Gamma \cap \mathcal{G}$ and $k \in \mathcal{K}$
where $\Gamma$ is the measurement site defined in \eqref{gamma def}, $\delta
= 10\%$ and $\mathrm{rand}$ is the function taking uniformly distributed
random numbers in $[-1, 1].$ The truncation number $N$ is $7$, which is
chosen by a trial-error process.

\subsubsection{The first approximation of the function $V$}

\label{subsec initial}

The first step of our method is to compute a vector valued function that
satisfies 
\begin{equation}
\partial_{\nu} V_0|_{\Gamma} = g_0 \quad \mbox{and} \quad V_0|_{\partial
\Omega} = g_1.  \label{boundary V0}
\end{equation}
This vector valued function $V_0$ is used in the change of variable $U = V -
V_0$ as in \eqref{change}. Moreover, to guarantee the fast convergence, we
will find $V_0$ such that it is close to the solution $V$. We call this
function $V_0$ the initial solution.

Since our target is to solve the nonlinear system \eqref{quasi-system}, it is natural to find $V_0 = (v^0_1, v^0_2, \dots, v^0_N)$ as the solution to a linear
system obtained by removing from \eqref{quasi-system} the nonlinear term, which is
\begin{equation}
    \left\{
        \begin{array}{rcll}
             \displaystyle\sum_{i = 1}^N s_{li} \Delta v_i^0(\mathbf{x}) 
%              + \sum_{i, j = 1}^N
% a_{lij}\nabla v_i(\mathbf{x}) \cdot \nabla v_j^0(\mathbf{x})
% \\
% \hspace{3cm}
 \displaystyle 
 +\sum_{i = 1}^N  
B_{li}(\mathbf{x}) \cdot \nabla v_i^0(\mathbf{x}) 
&=& 0,
             &\quad \x \in \Omega,\\
             V_0(\x)  &=& g_1(\x)  &\quad \x \in \partial \Omega,
             \\
             \partial_{\nu} V_0(\x)   &=& g_0(\x) &\quad \x \in \Gamma.
        \end{array}
    \right.
    \label{linear-system}
\end{equation}
 Since %
\eqref{linear-system}  is a system of linear partial
differential equations, we can solve it directly by the quasi-reversibility
method involving a Carleman weight function in the finite difference scheme.
That means, we minimize the following  functional 
\begin{equation}
    W \mapsto \int_{\Omega}  \mu_\lambda^2 \Big|\sum_{i = 1}^N s_{li} \Delta w(\mathbf{x}) +\sum_{i = 1}^N  
B_{li}(\mathbf{x}) \cdot \nabla W(\mathbf{x}) \Big|^2 + \epsilon\|W\|_{H^2(\Omega)}^2
\label{dropnonlinear}
\end{equation}
where $W = (w_1, w_2, \dots, w_N)$ is subject to the boundary conditions $W|_{\partial \Omega} = g_1$ and $\partial_{\nu }W|_{\Gamma} = g_0$. 
In \eqref{dropnonlinear} and also in this section, the used Carleman
weight function is $\mu_\lambda = e^{-\lambda(R + r)^2}e^{\lambda(z -
r)^2}$ where $\lambda = 1.1$ and $r = 1.5$ and the regularization parameter $\epsilon = 10^{-6}$. Even though in theory, the value
of $\lambda$ is large. However, we have discovered computationally that a
reasonable value $\lambda = 1.1$ works well. So, we use this $\lambda$.
These observations coincide with those of our previous works on the
numerical studies of the convexification \cite%
{KlibanovLiZhang:ip2019,KhoaKlibanovLoc:SIAMImaging2020}.
This Carleman weight function is used for all numerical tests in this section.

We refer the reader to \cite{LeNguyen:2020,Nguyen:CAMWA2020} for details in the
implementation of the quasi-reversibility method to solve a system of linear
partial differential equations with Cauchy boundary data. 

\subsubsection{The minimizing sequence}

For the simplification in implementation, we skip the step of changing the
variable $U = V - V_0$ as in \eqref{change}. Let $U = V - V_0$ where $V_0 =
(v_1^0, \dots, v_N^0)$ is the vector valued function found in section \ref%
{subsec initial}. Then, due to \eqref{6.11}, we set the cost functional as 
%\begin{equation}
%	\sum_{i = 1}^N s_{li} \Delta u_i(\x) = \sum_{i = 1}^N s_{li} \Delta v_i^0 - \sum_{i, j = 1}^N a_{lij} \nabla (u_i + v_i^0) \cdot \nabla (u_j + v_j^0) - \sum_{i = 1}^N B_{li}\nabla (u_i + v_i^0)
%\end{equation}
%for all $\x \in \Omega$, $l \in \{1, \dots, N\}$.
%Hence, to compute $U$, by the convexification method, we minimize the functional
\begin{equation*}
J(U) = \sum_{l = 1}^N \int_{\Omega} \mu_\lambda^2\Big|\sum_{i = 1}^Ns_{li}
\Delta v_i - \sum_{i, j = 1}^N a_{lij} \nabla v_i \cdot \nabla v_j - \sum_{i
= 1}^N B_{li}\nabla v_i\Big)\Big|^2 d\mathbf{x} + \epsilon
\|V\|^2_{H^2(\Omega)^N}.
\end{equation*}
The finite difference version of $J$ is 
\begin{multline}
J(V) = h^3\sum_{\mathfrak{i, j, l = 1}}^{N_\mathbf{x}}\sum_{l = 1}^N
\mu_\lambda^2(x_{\mathfrak{i}}, y_{\mathfrak{j}}, z_{\mathfrak{l}})\Big|%
\sum_{i = 1}^Ns_{li} \Delta v_i(x_{\mathfrak{i}}, y_{\mathfrak{j}}, z_{%
\mathfrak{l}}) \\
- \sum_{i, j = 1}^N a_{lij} \nabla v_i(x_{\mathfrak{i}}, y_{\mathfrak{j}},
z_{\mathfrak{l}}) \cdot \nabla v_j(x_{\mathfrak{i}}, y_{\mathfrak{j}}, z_{%
\mathfrak{l}}) - \sum_{i = 1}^N B_{li}(x_{\mathfrak{i}}, y_{\mathfrak{j}},
z_{\mathfrak{l}})\nabla v_i(x_{\mathfrak{i}}, y_{\mathfrak{j}}, z_{\mathfrak{%
l}}) \Big|^2 d\mathbf{x} \\
+ \epsilon h^3 \sum_{\mathfrak{i, j, l = 1}}^{N_\mathbf{x}} \sum_{l = 1}^N %
\Big(|v_l(x_{\mathfrak{i}}, y_{\mathfrak{j}}, z_{\mathfrak{l}})|^2 + |\nabla
v_l(x_{\mathfrak{i}}, y_{\mathfrak{j}}, z_{\mathfrak{l}})|^2 + |\Delta
v_l(x_{\mathfrak{i}}, y_{\mathfrak{j}}, z_{\mathfrak{l}})|^2 \Big).
\label{8.4}
\end{multline}

\begin{Remark}
In our computation, $\epsilon = 10^{-6}$ for all tests. Also, in \eqref{8.4}%
, the regularity term is set to be $\epsilon\|U\|_{H^2(\Omega)^N}$ instead
of $H^p(\Omega)^N$. We observe numerically that using the norm $%
\|U\|_{H^2(\Omega)^N}$ already provides satisfactory numerical solutions.
So, it is not necessary for us to choose norm in $H^p(\Omega)^N$. This
observation significantly reduces the expensive efforts in implementation.
\end{Remark}

We now mention that to speed up the minimization procedure, we need to
compute the gradient $DJ_1$ of the discrete functional $J_1$ in \eqref{8.4}
above. Having the expression for the gradient via an explicit formula
significantly reduces the computational time. We have derived such a formula
using the technique of Kronecker deltas, which has been outlined in \cite{Kuzhuget:AA2010}. For brevity we do not provide this formula here.

\subsubsection{Numerical examples}

We perform three (3) tests.

\textit{Test 1.} We first consider the case of detecting one target with
high dielectric constant. The dielectric constant of the medium is given by 
\begin{equation*}
c_{\mathrm{true}} =  \left\{  
\begin{array}{ll}
5 & \mbox{if } 0.6x^2 + y^2 + (z + 0.7)^2 \leq 0.2^2, \\ 
1 & \mbox{otherwise}.%
\end{array}
\right. 
\end{equation*}

\begin{figure}[h!]
\begin{center}
\subfloat[The true 3D
image]{\includegraphics[width=.3\textwidth]{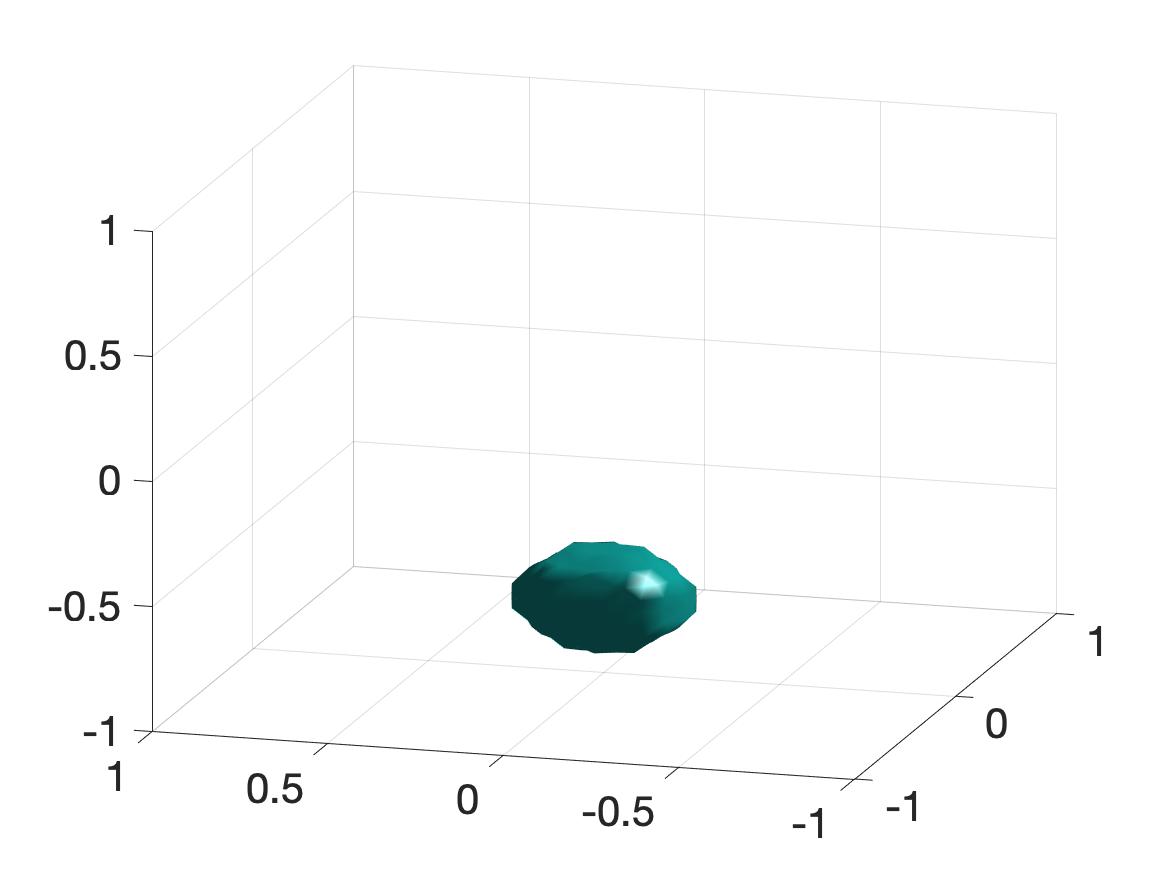}}  \quad \quad  %
\subfloat[The reconstructed 3D
image]{\includegraphics[width=.3\textwidth]{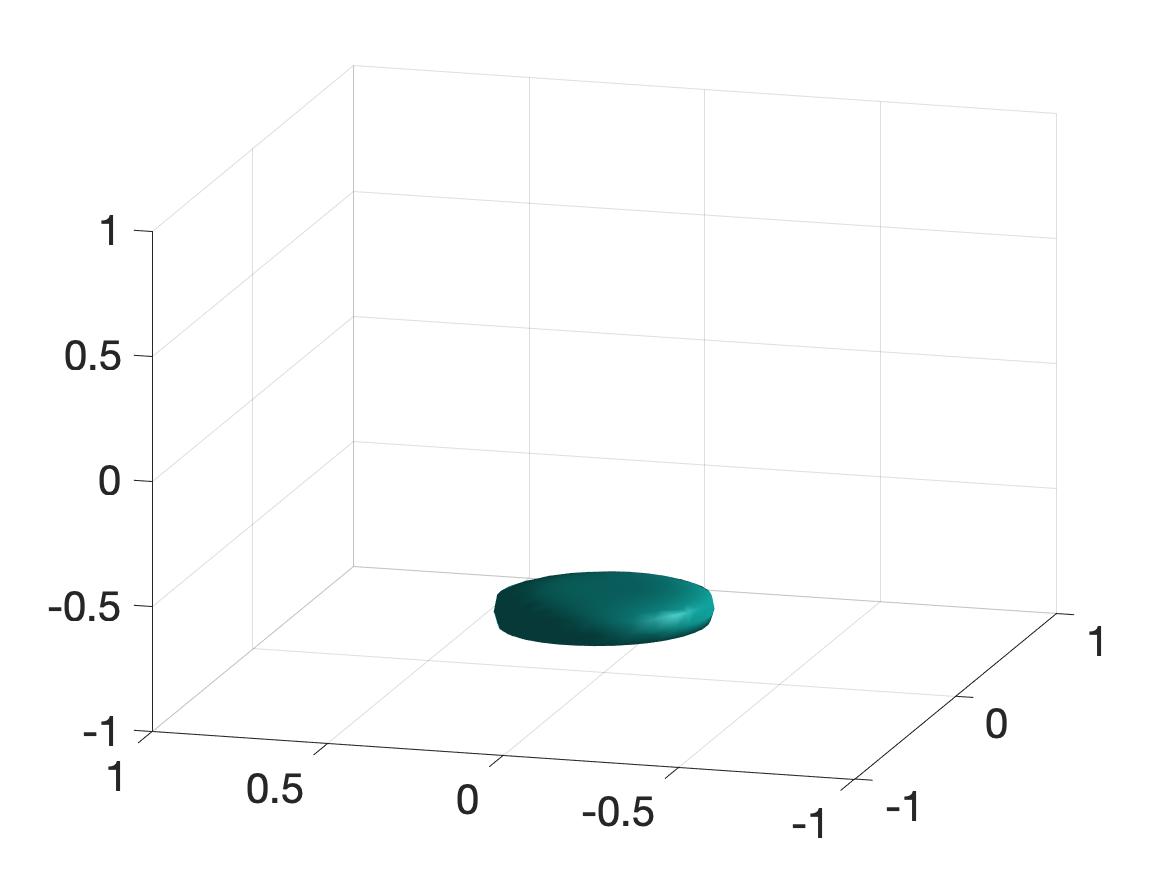}}
\par
\subfloat[The cross sections $\{z = -0.7\}$ and $\{y = 0\}$ of $c_{\rm
true}$]{\includegraphics[width=.35\textwidth]{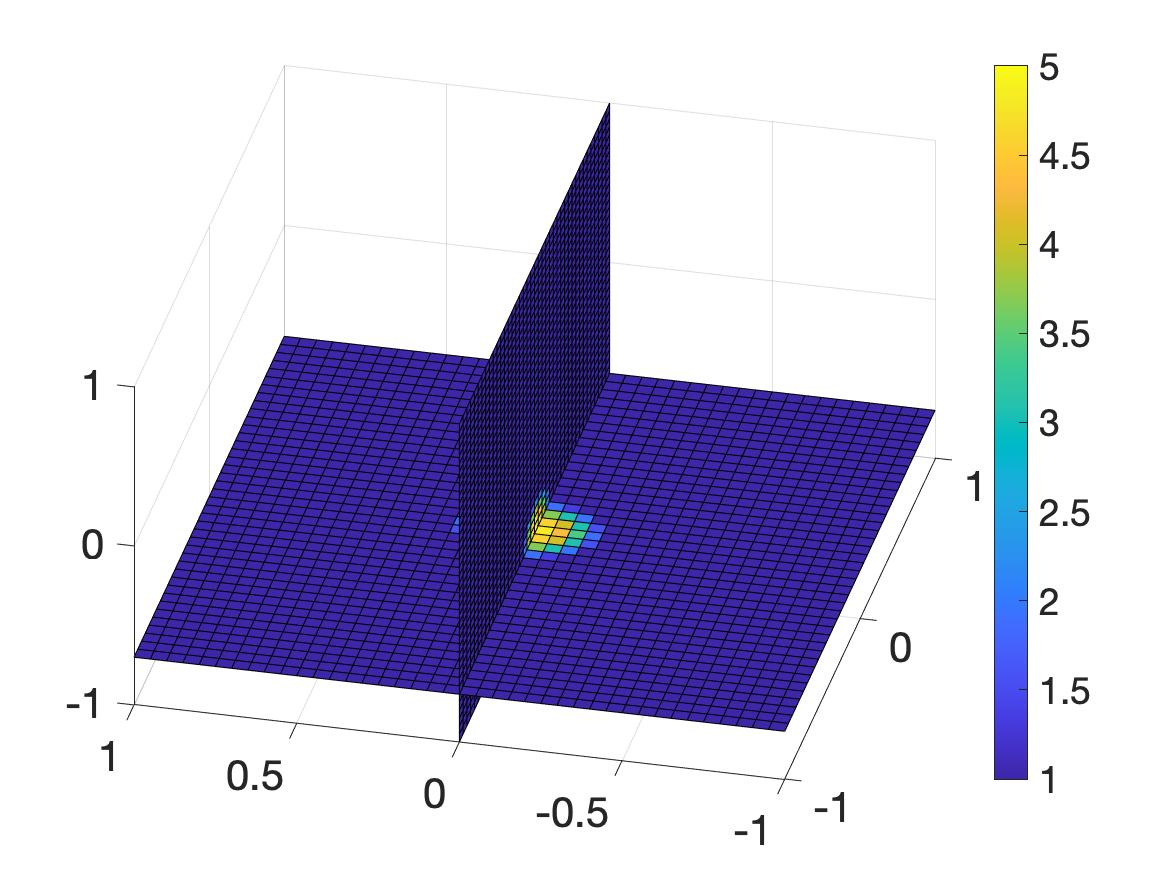}}  \quad \quad  %
\subfloat[The cross sections $\{z = -0.7\}$ and $\{y = 0\}$ of $c_{\rm
comp}$]{\includegraphics[width=.35\textwidth]{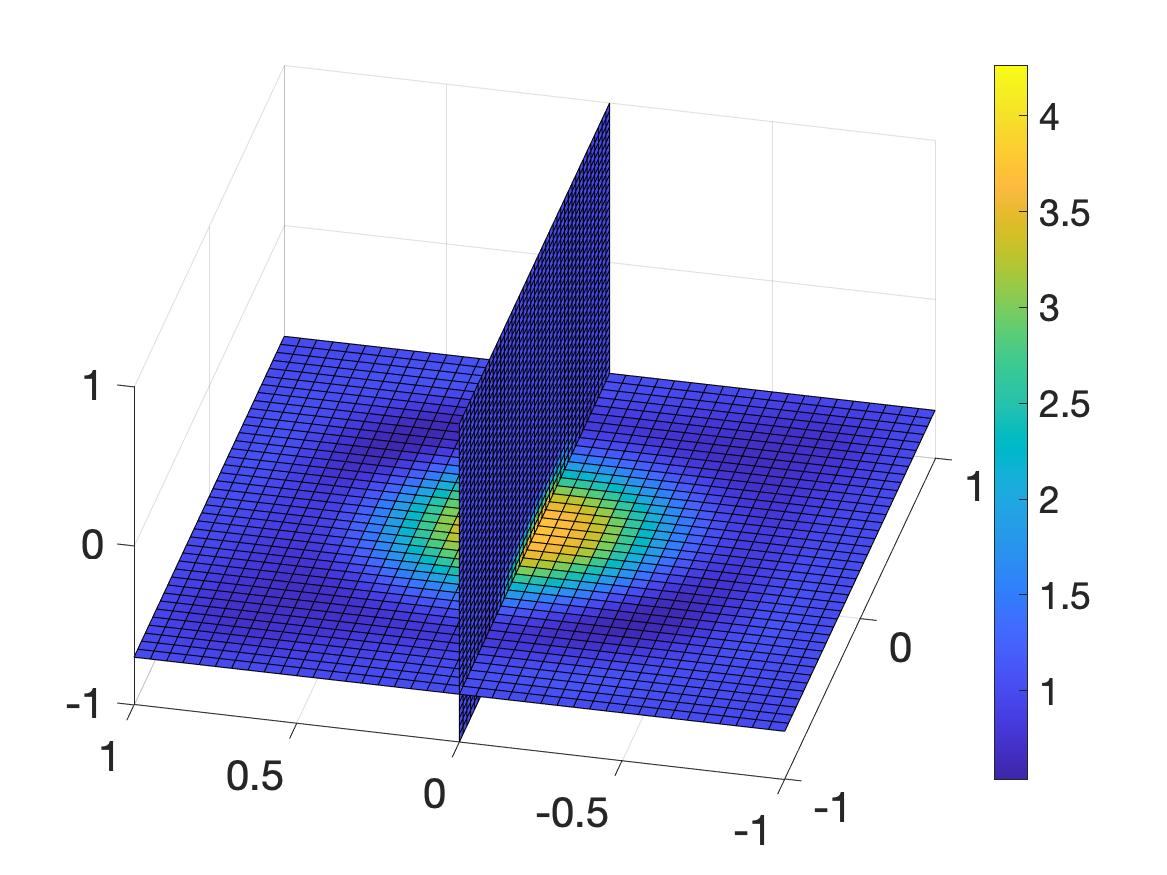}}  
\end{center}
\caption{Test 1. The function $c_{\mathrm{true}}$ and its reconstruction $c_{%
\mathrm{comp}}$ from noisy data with noise level of $10\%.$ }
\label{fig test 1}
\end{figure}

The true and computed dielectric constants are displayed in Figure \ref{fig
test 1}. It is obvious that the location of the target is detected
accurately. The reconstructed shape is somewhat acceptable. The computed
maximal value of the dielectric constant is $4.26$ (relative error 14.8\%)).

\textit{Test 2.} We test our method when the true dielectric constant is
given by 
\begin{equation*}
c_{\mathrm{true}} =  \left\{  
\begin{array}{ll}
3 & \mbox{if } 0.35^2 \leq x^2 + y^2 \leq 0.5^2 \mbox{ and } -0.8 \leq z
\leq -0.65 \\ 
1 & \mbox{otherwise}.%
\end{array}
\right. 
\end{equation*}
The shape of the dielectric constant is a ring.

\begin{figure}[h!]
\begin{center}
\subfloat[The true 3D image of the
ring]{\includegraphics[width=.3\textwidth]{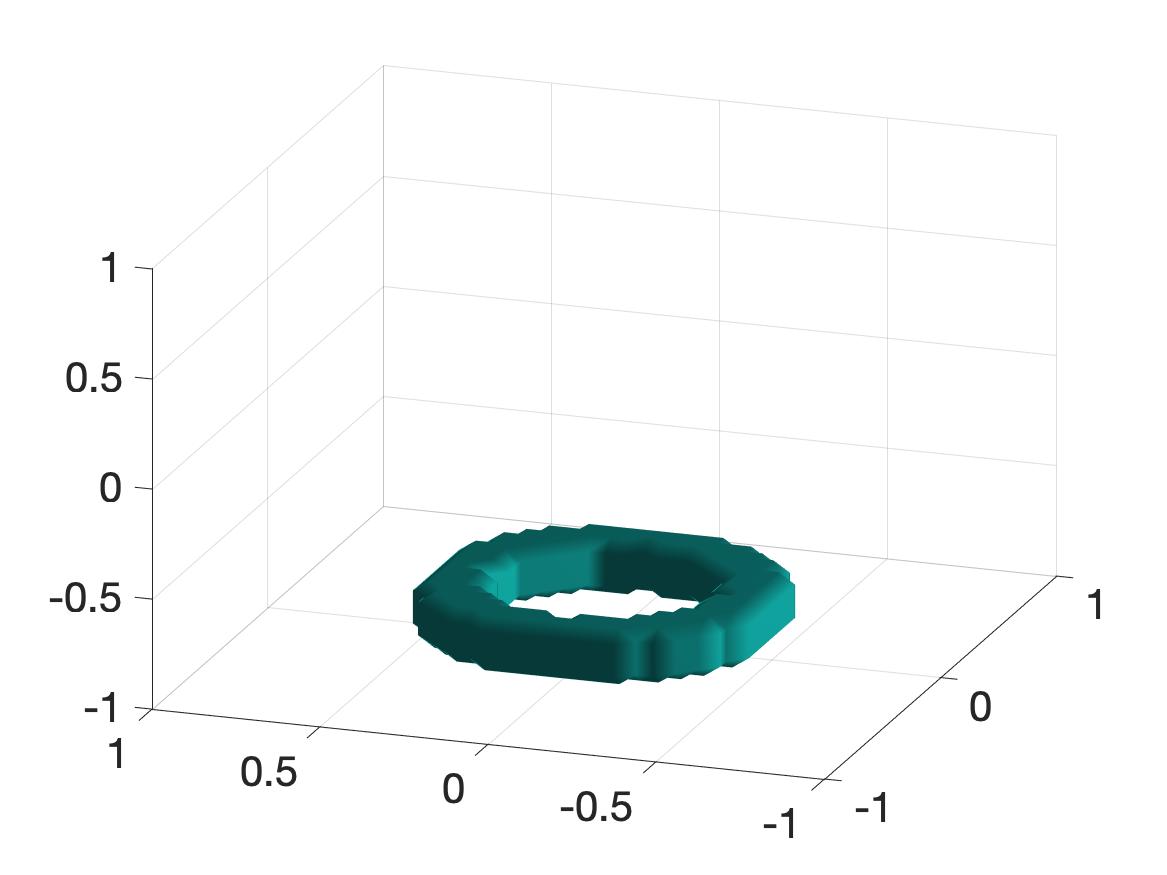}}  \quad \quad  %
\subfloat[The reconstructed 3D image of the
ring]{\includegraphics[width=.3\textwidth]{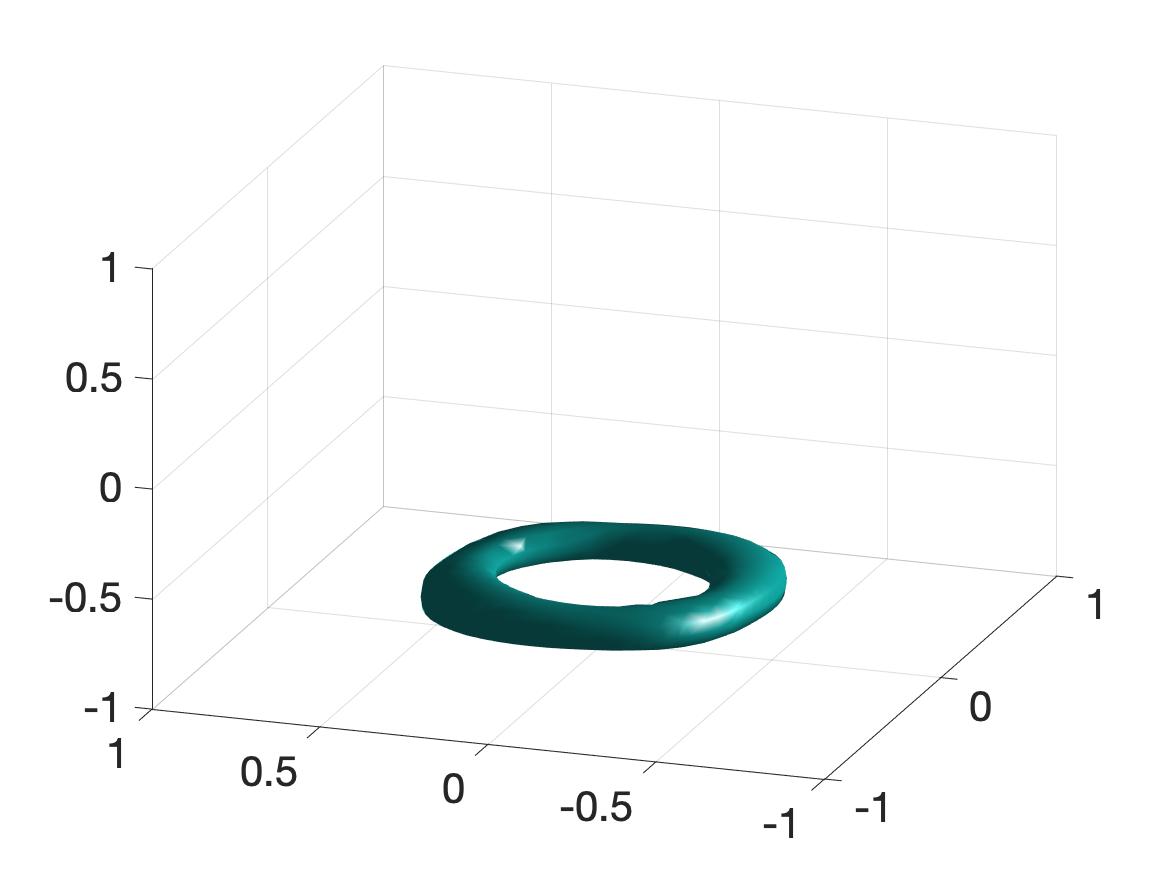}}
\par
\subfloat[The cross sections $\{z = -0.7\}$ and $\{y = 0\}$ of $c_{\rm
true}$]{\includegraphics[width=.35\textwidth]{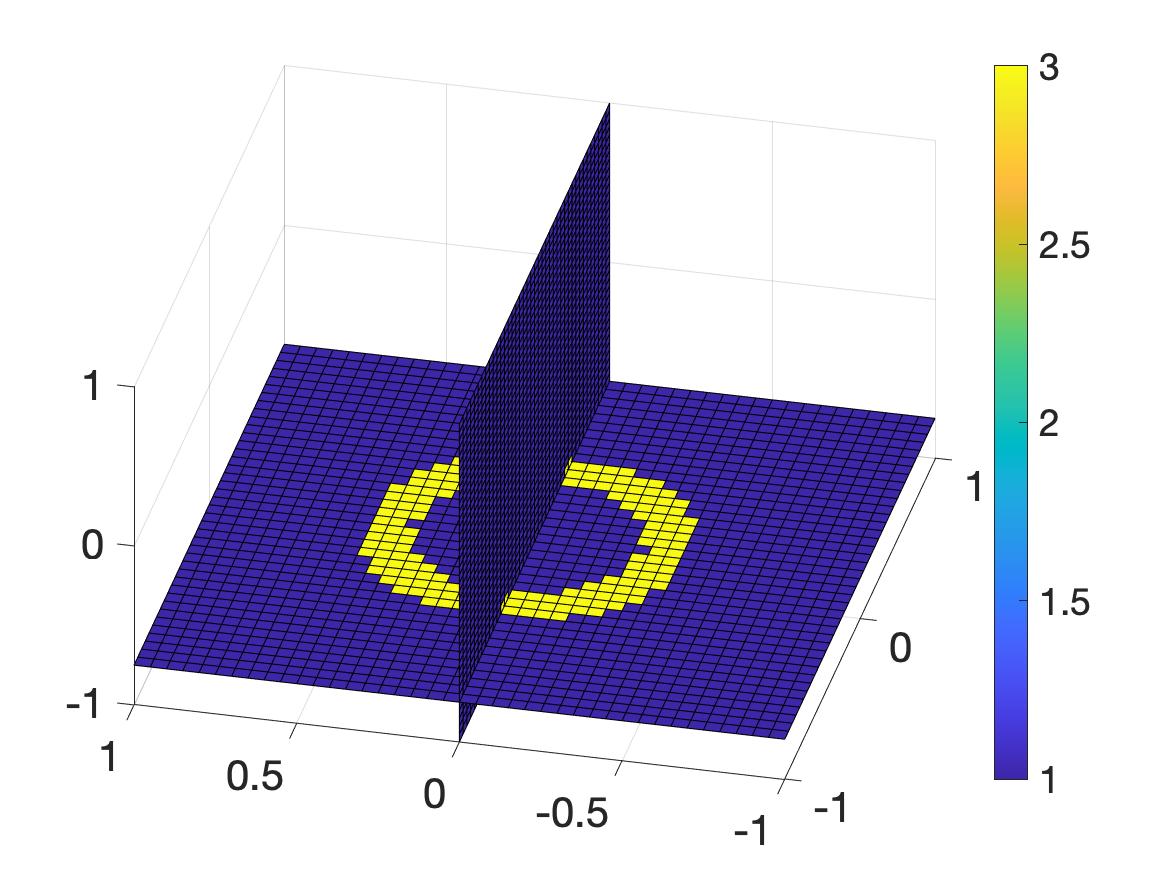}}  \quad \quad  %
\subfloat[The cross sections $\{z = -0.7\}$ and $\{y = 0\}$ of $c_{\rm
comp}$]{\includegraphics[width=.35\textwidth]{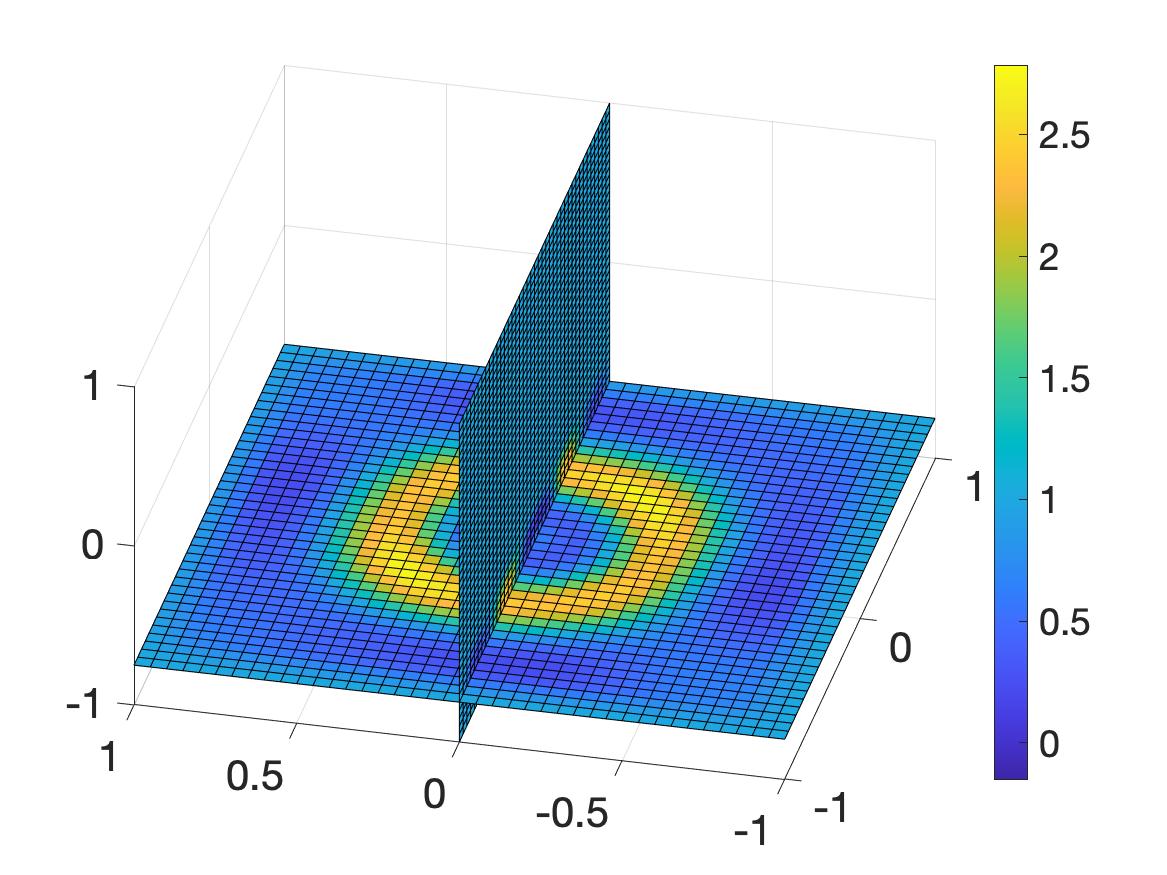}}  
\end{center}
\caption{Test 2. The function $c_{\mathrm{true}}$ and its reconstruction $c_{%
\mathrm{comp}}$ from noisy data with noise level of $10\%.$ }
\label{fig test 2}
\end{figure}

The true and computed dielectric constants are displayed in Figure \ref{fig
test 2}. It is an evident that the dielectric constant is computed
successfully. The ``ring" shape is clearly detected. The computed maximal
value of the dielectric constant is $2.7809$ (relative error 7.3\%)).

We now consider the direct optimization without using the Carleman weight
function. That means we apply the same procedure to compute the dielectric
constant except taking $\lambda = 0.$ The numerical result in Figure \ref%
{fig nolambda} show that without the Carleman weight function involving, we
reconstruction is poor.

\begin{figure}[h!]
\begin{center}
\subfloat[The 3D image computed without using the convexification
method]{\includegraphics[width=.3\textwidth]{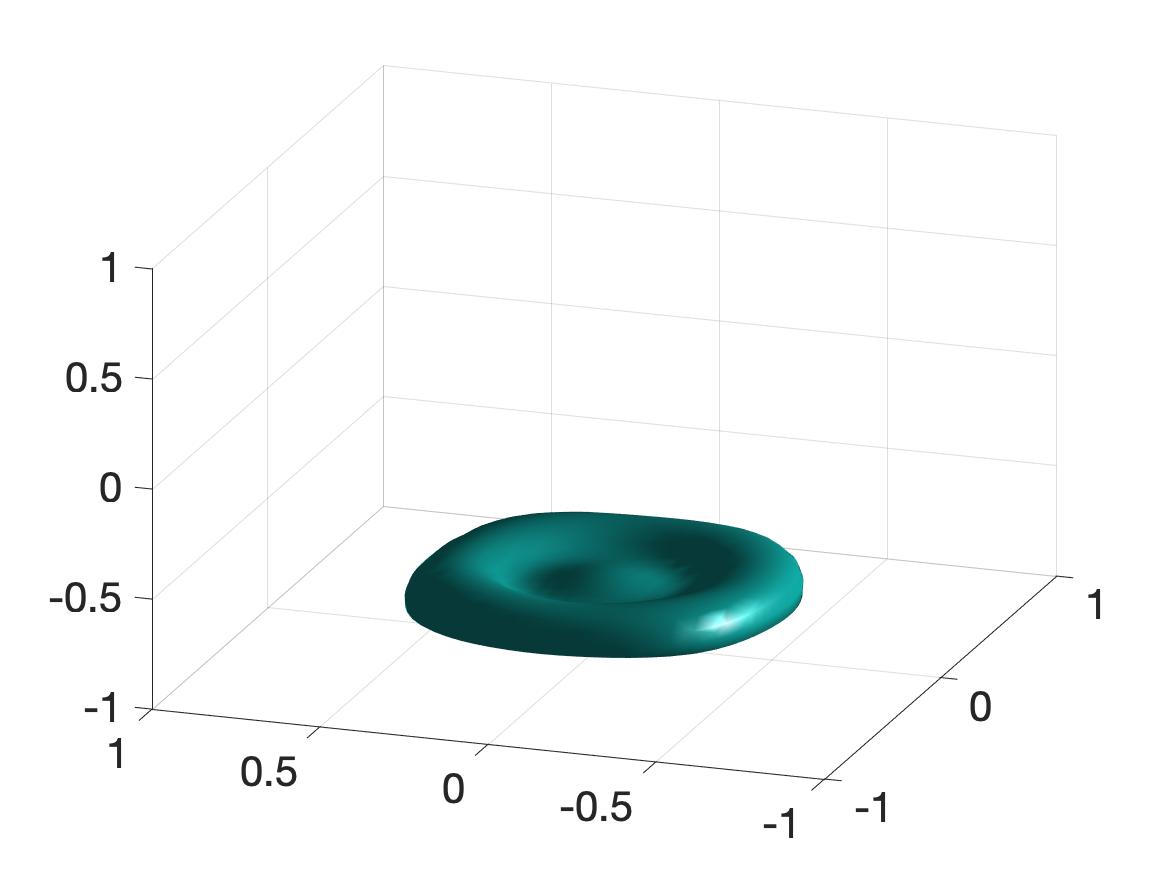}}  \quad
\quad  \subfloat[The cross sections $\{z = -0.75\}$ and $\{y = 0\}$ of
$c_{\rm comp}$]{\includegraphics[width=.3\textwidth]{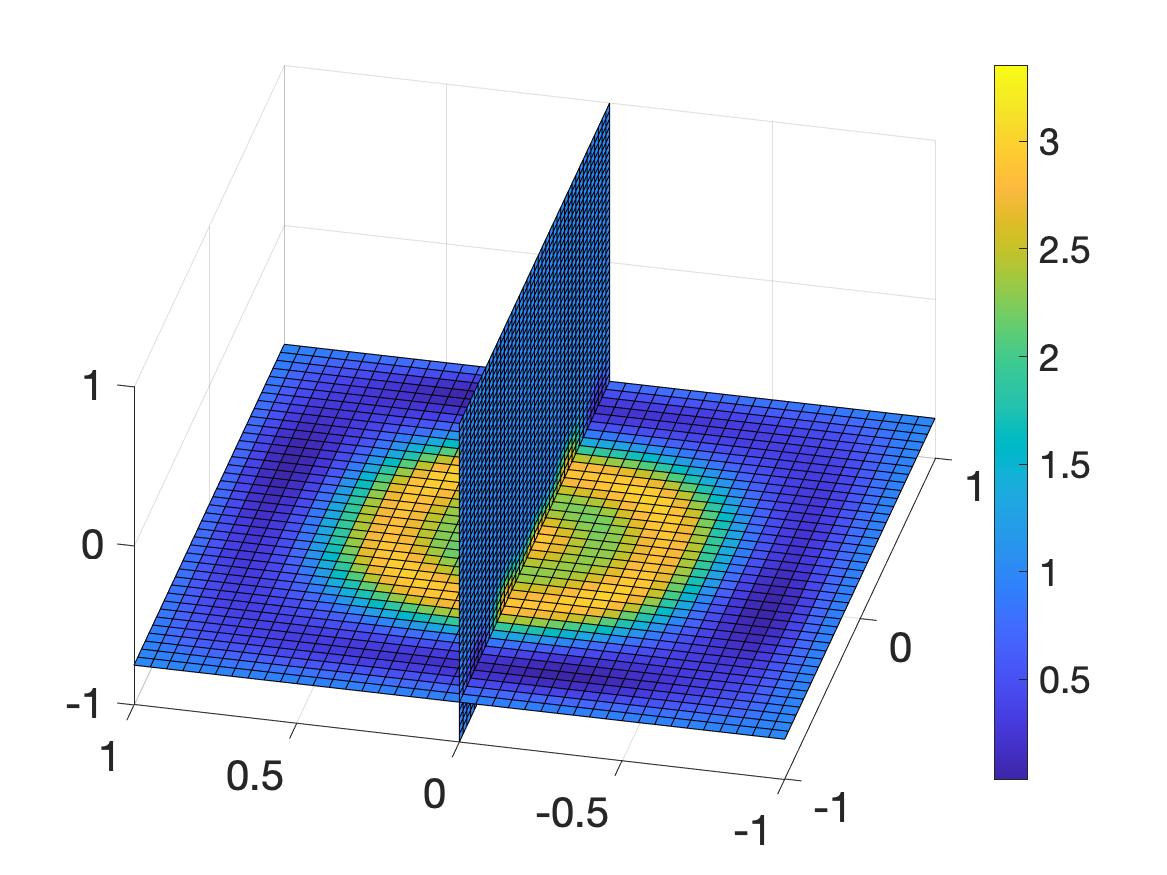}} 
\end{center}
\caption{Test 2. The function $c_{\mathrm{true}}$ and its reconstruction $c_{%
\mathrm{comp}}$ from noisy data with noise level of $10\%$ without using
Carleman weight function. It is evident that in this case, the ``ring" shape
cannot be reconstructed well.}
\label{fig nolambda}
\end{figure}

\textit{Test 3.} We consider dielectric constant with a more complicate
geometry. The graph of the dielectric constant is a letter $Y$ located on
the plane $z = -.7$ 
\begin{figure}[h!]
\begin{center}
\subfloat[The true 3D image of the letter
$Y$]{\includegraphics[width=.3\textwidth]{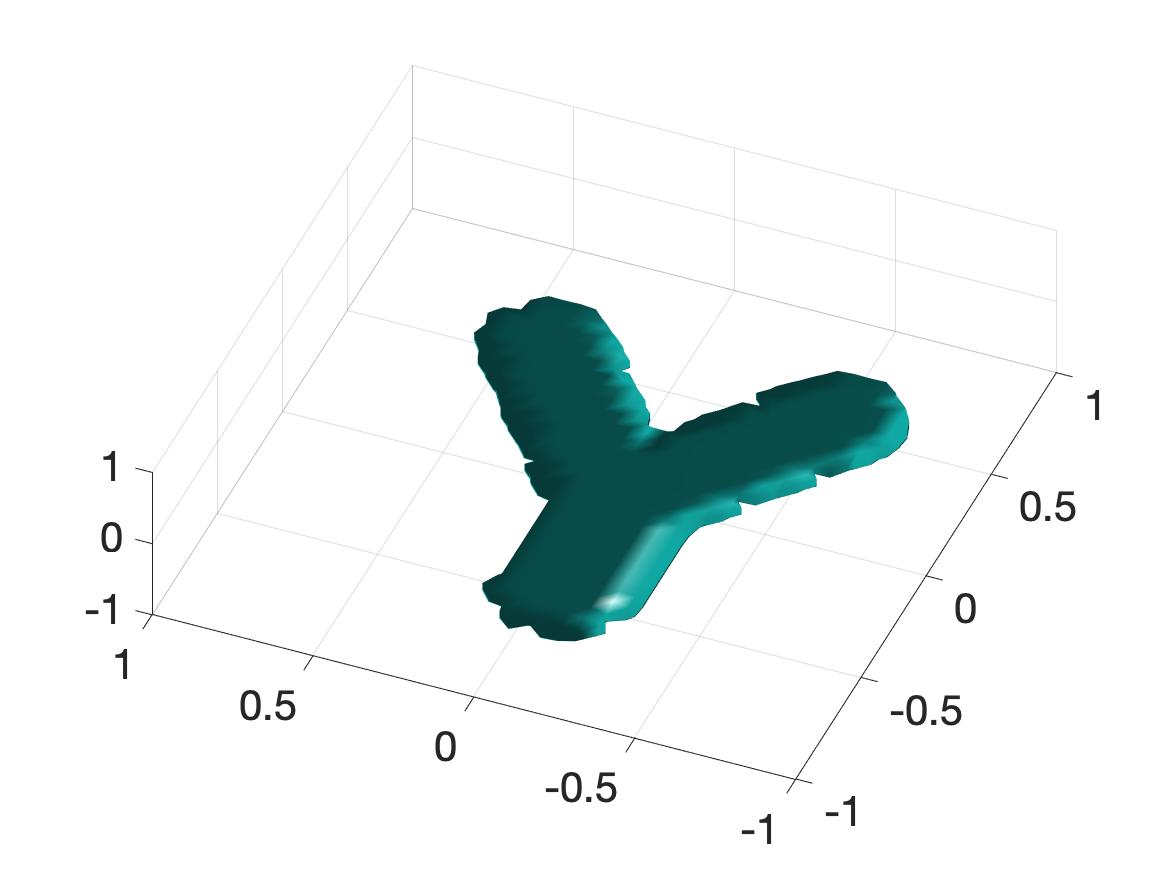}}  \quad \quad  %
\subfloat[The reconstructed 3D image of the letter
$Y$]{\includegraphics[width=.3\textwidth]{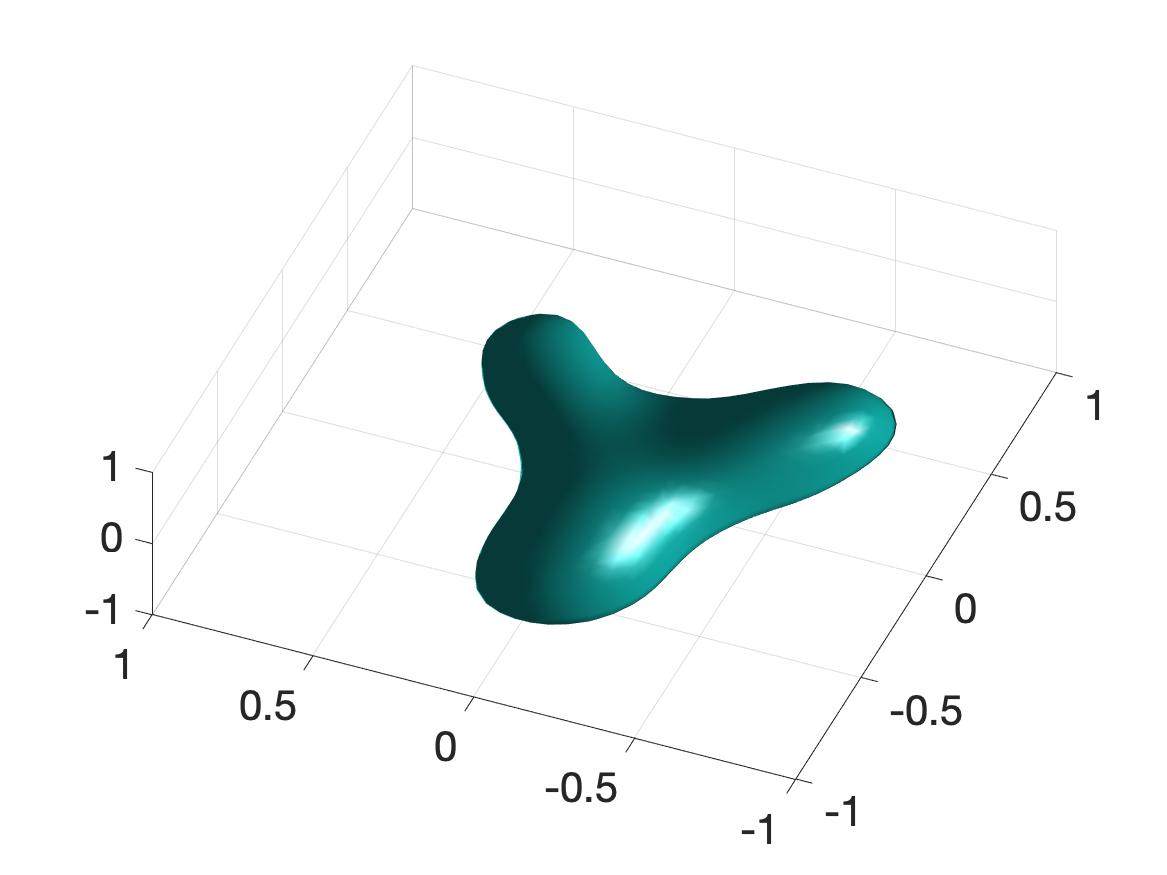}}
\par
\subfloat[The cross sections $\{z = -0.7\}$ and $\{y = 0\}$ of $c_{\rm
true}$]{\includegraphics[width=.35\textwidth]{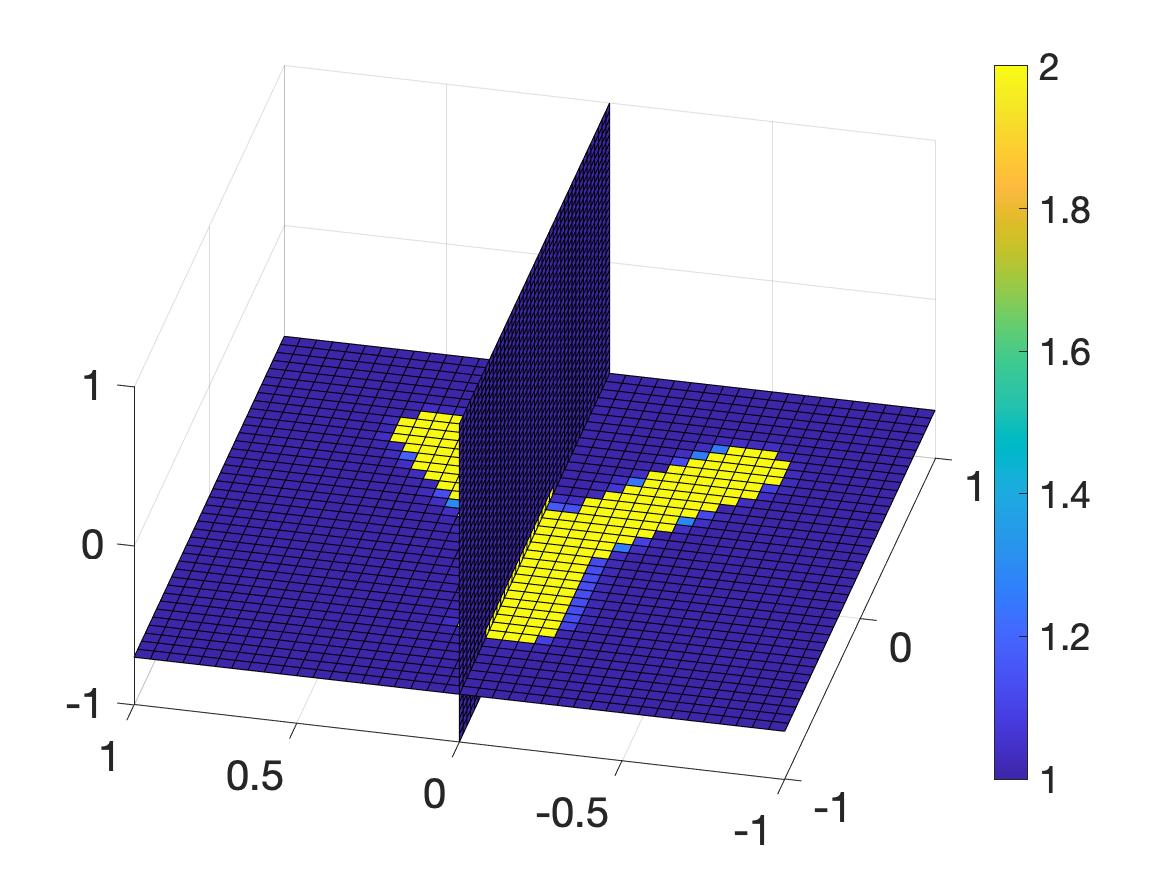}}  \quad \quad  %
\subfloat[The cross sections $\{z = -0.7\}$ and $\{y = 0\}$ of $c_{\rm
comp}$]{\includegraphics[width=.35\textwidth]{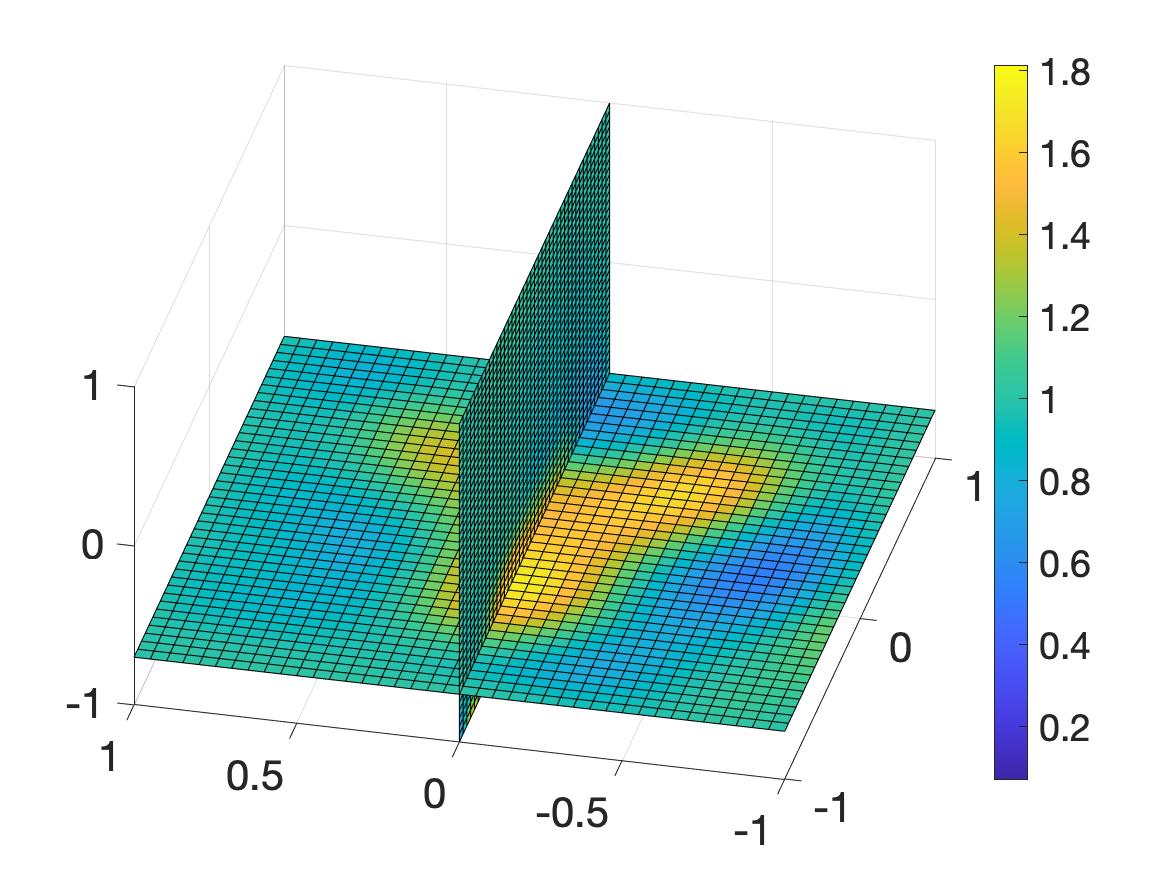}}  
\end{center}
\caption{Test 3. The function $c_{\mathrm{true}}$ and its reconstruction $c_{%
\mathrm{comp}}$ from noisy data with noise level of $10\%.$ }
\label{fig test 3}
\end{figure}
The true and constructed dielectric constants are displayed in Figure \ref%
{fig test 3}. We observe that our method can detect the shape of the letter $%
Y$ clearly. Moreover, the reconstructed value of the dielectric constant is
acceptable. The computed maximal value of the function $c$ is $1.8116$
(relative error 9.4\%).

\section{Concluding remarks}

\label{sec 9}

In the first part of this paper, we proved the convergence of the gradient
descent method to find the minimizer of a functional which is
strictly convex on a ball in a Hilbert space, rather than on the whole space.
 This is a new result, compared with previously obtained ones by
our research team for the case of a more complicated gradient projection
method. Then we used the convexification method and gradient descent method
to solve a boundary value problem of quasi-linear PDE with both Dirichlet
and Neumann data. We proved that this approach provides good numerical
solutions as the noise tends to zero. In the second part of the paper, we
applied the theoretical results of the first part to solve an inverse
scattering problem. To solve this inverse problem, we derive an approximate
mathematical model, which is the Cauchy problem for a coupled system of
quasilinear elliptic partial differential equations. Then, we apply the
convexification and the gradient descent method to solve this system. 
Numerical results for the inverse scattering problem demonstrate a
good reconstruction quality.

\section*{Acknowledgments}

 The authors sincerely thank Dr. Michael
V. Klibanov for many fruitful discussions.
The work is supported in part by US Army Research Laboratory and US Army Research
Office grant W911NF-19-1-0044 and by funds provided by the Faculty Research Grant program at UNC Charlotte, Fund No. 111272.

%\bibliographystyle{plain}
%%\bibliography{../../../../../../mybib}
%\bibliography{mybib}

%\section*{Statements and Declarations}
%
%\noindent{\bf Funding.} The work is supported in part by US Army Research Laboratory and US Army Research
%Office grant W911NF-19-1-0044 and by funds provided by the Faculty Research Grant program at UNC Charlotte, Fund No. 111272.
%
%\noindent{\bf Competing interests.} The authors have no relevant financial or non-financial interests to disclose.
%
%\noindent{\bf Author Contributions.} All authors contributed to the manuscript equally.
%
%\noindent{\bf Data Availability.} The data is generated computationally by solving a number of Partial Differential Equations. It is available from the corresponding author on reasonable request.

\end{document}